\documentclass{pspum-l}

\usepackage[dvips]{graphicx}
\usepackage[latin1]{inputenc}
\usepackage[T1]{fontenc}
\usepackage{amssymb}

\newcommand{\abs}[1]{\ensuremath{\left\vert #1 \right\vert}}
\newcommand{\hdim}{\text{\textup{dim}}_{\text{\textup{H}}}}

\makeatletter
\renewcommand{\subsubsection}{\@startsection
  {subsubsection}%
  {3}%
  {0pt}%
  {-.5\baselineskip}%
  {-\fontdimen2\font
    plus -\fontdimen3\font}%
  {\normalfont\itshape}%
}
\makeatother

\DeclareFontFamily{OMS}{rsfs}{\skewchar\font'177}
\DeclareFontShape{OMS}{rsfs}{m}{n}{%
  <5> rsfs5
  <6> <7> rsfs7
  <8> <9> <10> rsfs10
  <10.95> <12> <14.4> <17.28> <20.74> <24.88> rsfs10
}{}
\DeclareSymbolFont{rsfscript}{OMS}{rsfs}{m}{n}
\DeclareSymbolFontAlphabet{\mathrsfs}{rsfscript}
\newcommand{\scr}[1]{\mathrsfs{#1}}  
\DeclareMathOperator{\di}{diam}
\DeclareMathOperator{\SL}{SL}
\DeclareMathOperator{\codim}{codim}

\copyrightinfo{2003}{American Mathematical Society}

\newtheorem{theorem}{Theorem}[section]
\newtheorem{lemma}[theorem]{Lemma}

\newtheorem{cor}[theorem]{Corollary}

\theoremstyle{definition}

\theoremstyle{remark}

\numberwithin{equation}{section}

\begin{document}

\title{Hausdorff Dimension and Diophantine Approximation}

\author{M. Maurice Dodson}
\address{University of York\\
  Heslington \\
  York \\
  YO10 5DD \\
  UK
}
\email{mmd1@york.ac.uk}

\author{Simon Kristensen}
\address{University of York\\
  Heslington \\
  York \\
  YO10 5DD \\
  UK
}
\email{sk17@york.ac.uk}
\thanks{Research funded by EPSRC grant no. GR/N02832/01 with
  additional support from INTAS grant no. 001--429.}




\maketitle

\section{Introduction }

Dimension provides an indication of the size and complexity of a set
and various kinds, such as box counting, packing and Hausdorff
dimensions, play an important role in the study of
fractals~\cite{FalcFG}.  For example, the Hausdorff dimension of the
Cantor middle third set is $\log 2/\log 3$ (proved by Hausdorff in his
seminal paper~\cite{Hausdorff19}), that of the Koch snowflake curve is
$\log 4/\log 3$ and it has recently been shown that the boundary of
the Mandelbrot set, a very complicated set of Lebesgue measure 0 or
\emph{null} set in the complex plane with topological dimension 1, has
Hausdorff dimension 2~\cite{shishikura1998}.  On the other hand,
Diophantine approximation is a quantitative analysis of rational
approximation and so, at least at first sight, is less geometrical.
The purpose of this article is to show that Hausdorff dimension plays
an important part in this theory too.

In order to keep the article accessible, the emphasis is on
approximation of real numbers by rationals and the less well known
topic of approximation of complex numbers by ratios of Gaussian
integers. The more general theory, which recently has seen some
spectacular advances, will be referred to and some applications
sketched. The article is organised as follows.  We begin with a brief
treatment of Hausdorff measure and Hausdorff dimension. We then
explain some of the principal results in Diophantine approximation and
the Hausdorff dimension of related sets, originating in the pioneering
work of Vojt{\v{e}}ch Jarník~\cite{NSJobit}.  We conclude with some
applications of these results to the metrical structure of exceptional
sets associated with some famous problems. It is not intended that all
the recent developments be covered but they can be found in the
references cited.

\section{Hausdorff measure and dimension}

Felix Hausdorff introduced the notion of Hausdorff dimension in a
remarkable and influential paper~\cite{Hausdorff19} that extended
Carath\'eodory's approach to Lebesgue measure~\cite{Carat14} in a
simple but far-reaching way.  (Kahane's Foreword to the book
\emph{Fractals} ~\cite{Cherbit91} includes a short and moving
biography of Hausdorff.)  Dimension had been taken to be a
non-negative integer but by a simple observation, which Hausdorff
described modestly as `a small contribution', he modified
Cara\-th\'eodory's definition of measure to obtain a measure
associated with a dimension that could be any non-negative real
number.  We shall assume a knowledge of Lebesgue measure and as usual,
we shall often say \emph{almost no} to indicate a null set -- thus
almost no numbers are rational -- and we shall say \emph{almost all}
to indicate a set whose complement is null, so that almost all numbers
are irrational.

For familiar sets such as the interval, circle and the plane, the
Hausdorff dimension\ (defined below) coincides with the usual notion
of dimension and is respectively 1,\,1 and 2.  However, a significant
difference is that \emph{any\/} set in Euclidean space has a Hausdorff
dimension (a non-measurable set in $\mathbb{R}^n$ has full Hausdorff
dimension $n$).  In particular, null sets, such as Cantor's middle
third set or the set of badly approximable numbers
(see~\S\ref{sec:toa}), have a Hausdorff dimension and this gives a way
of discriminating between them.  It is also natural to study the
Hausdorff dimension of \emph{exceptional sets} which are sets
associated with the invalidity of some result, making it desirable
that they be null.  A brief and more or less self-contained account of
Hausdorff measure and dimension is now given (more detailed
expositions can be found
in~\cite{MDAM,FalcGFSshort,FalcFG,Fed,MattilaGS,Rogers}).

\subsection{Hausdorff measure} 

Carathéodory's approach to the measure of a set $E$ in $\mathbb{R}^n$
was based on `approximating' $E$ by countable covers consisting of
small `simple' sets $U$ in $ \mathbb{R}^n$.  Hausdorff's idea was to
introduce for a given cover, $\scr{C}$ say, of $E$ the sum (sometimes
termed the \emph{$s$-length} of the cover $\scr{C}$)
\begin{equation*}
  \ell^s({\scr{C}}):=
  \sum_{U\in{\scr{C}}} (\di U)^{s},
\end{equation*}
where $\di U=\sup\{\vert \mathbf{x}-\mathbf{y}\vert_2\colon
\mathbf{x},\mathbf{y} \in U\}$ is the diameter of $U$ ($\vert
\mathbf{x}-\mathbf{y}\vert_2$ is the usual Euclidean distance between
$\mathbf{x}$ and $\mathbf{y}$) and where $s$ \emph{is a non-negative
  real number that is not necessarily an integer}.  Hausdorff
considered a monotonic function $l$ which allows more discrimination
but for simplicity we shall stick to the more familiar widely used
special case $l(\di U)=(\di U)^s$, associated with what is now usually
called the Hausdorff dimension but also sometimes called the
Hausdorff-Besicovitch dimension~\cite{Rogers}.  The possibly infinite
number $\ell^s({\scr{C}})$ gives an indication of the number of
subsets $U$ in $\scr{C}$ needed to cover $E$.  In order to effect the
approximation, the diameter of the sets $U$ in the cover is restricted
to be at most $\delta>0$.

Let 
\begin{equation*}
  \mathcal{H}^s_\delta(E):=\inf_{\scr{C}_\delta} 
  \sum_{U\in{\scr{C}}_\delta} (\di
  U)^{s} =\inf \ell^s({\scr{C_\delta}}),
\end{equation*}
where the infimum is taken over all covers $\scr{C}_\delta$ of $E$ by
sets $U$ with $\di U \leqslant \delta$; such covers are called
$\delta$-covers.  For a point $\mathbf{x}$,
$\mathcal{H}^s_\delta(\{\mathbf{x}\})=1$ when $s=0$ and vanishes when
$s>0$.  As $\delta$ decreases, $\mathcal{H}^s_\delta$ can only
increase as there are fewer $U$'s available, \emph{i.e.}, if
$0<\delta<\delta'$, then
\begin{equation*}
  \mathcal{H}^s_{\delta'}(E)\leqslant \mathcal{H}^s_{\delta}(E). 
\end{equation*}
The set function $\mathcal{H}^s_\delta$ is an outer measure on
$\mathbb{R}^n$ but the limit $\mathcal{H}^s$ (which can be infinite)
as $\delta \to 0$, given by
\begin{equation}
\label{eq:HHd}
  \mathcal{H}^s(E)=\lim_{\delta\to 0}\mathcal{H}^s_\delta (E) \ = 
  \sup_{\delta>0} \mathcal{H}^s_\delta (E)\in [0,\infty],
\end{equation}
is better behaved.  From its construction by covers,
$\mathcal{H}^s(E)\leqslant \mathcal{H}^s(F)$ for any $E\subset F$ and
indeed is subadditive and a regular outer measure.  The restriction to
the $\sigma$-field of $\mathcal{H}^s$ measurable sets (which includes
open and closed sets, limsup and liminf sets and $G_\delta$ and
$F_\sigma$ sets) is usually called the \emph{Hausdorff $s$-dimensional
measure}.  Hausdorff $1$-dimensional measure coincides with
1-dimensional Lebesgue measure and in higher dimensions, Hausdorff
$n$-dimensional measure is comparable to $n$-dimensional Lebesgue
measure, \emph{i.e.},
\begin{equation*}
  \mathcal{H}^n(E)\asymp \vert E\vert ,
\end{equation*}
where $\vert E\vert $ is the Lebesgue measure of $E$ and where for $a,b>0$,
$a\asymp b$ means there exist constants $c,c'>0$ such that $a\leqslant
c b\leqslant c' a$ or $a=O(b)$, $b=O(a)$ in Landau's $O$-notation.
Thus a set of positive $n$-dimensional Lebesgue measure has positive
Hausdorff $n$-measure.

Because it is defined in terms of the diameter of the covering sets,
Hausdorff $s$-measure is unchanged by restriction to closed, convex or
open sets.  It is also unchanged by isometries and so in particular by
translations and rotations.  It is, however, affected by scaling in
the natural way (as are fractals): for any $r\geqslant 0$,
\begin{equation*}
  \mathcal{H}^s(r E)=r^s \mathcal{H}^s(E).
\end{equation*}

\subsection{Hausdorff dimension} 
\label{sec:HD}
Zero-dimensional Hausdorff measure $\mathcal{H}^0(E)$ is simply
counting measure; thus the Hausdorff $s$-measure of a set of $k$
points is $k$ when $s=0$ and 0 for $s>0$.  This pattern is typical.
When the set $E$ is infinite, $\mathcal{H}^s(E)$ is either 0 or
$\infty$, except for possibly one value of $s$. To see this, the
definition of $\mathcal{H}^{s}_\delta(E)$ implies that there is a
$\delta$-cover $\scr{C}_\delta$ of $E$ such that
\begin{equation*}
  \sum_{C\in\scr{C}_\delta}(\di C)^{s} \leqslant
  \mathcal{H}^{s}_\delta(E)+1 \leqslant \mathcal{H}^{s}(E)+1 \leqslant
  \infty.
\end{equation*}
Suppose that $\mathcal{H}^{s_{_0}}(E)$ is finite and
$s=s_0+\varepsilon$, $\varepsilon>0$.  Then for each member $C$ of the
cover $\scr{C}_\delta$, $(\di C)^{s_0+\varepsilon}\leqslant
\delta^\varepsilon (\di C)^{s_0}$, so that the sum
\begin{equation*}
  \sum_{C\in\scr{C}_\delta}( \di C)^{s_0+\varepsilon}\leqslant
  \delta^\varepsilon \sum_{C\in\scr{C}_\delta}(\di C)^{s_0}.
\end{equation*}
Hence 
\begin{equation*}
  \mathcal{H}_\delta^{s_0+\varepsilon}(E)\leqslant
  \sum_{C\in\scr{C}_\delta} ( \di C)^{s_0+\varepsilon} \leqslant
  \delta^\varepsilon \sum_{C\in\scr{C}_\delta} (\di
  C)^{s_0}\leqslant\delta^\varepsilon( \mathcal{H}^{s_{_0}}(E)+1)
\end{equation*}
and so
\begin{equation*}
  0\leqslant \mathcal{H}^s(E)= \mathcal{H}^{s_{_0}+\varepsilon}(E)=
  \lim_{\delta\to0}\mathcal{H}_\delta^{s_{_0}+\varepsilon}(E)\leqslant
  \lim_{\delta\to0} \delta^\varepsilon( \mathcal{H}^{s_0}(E)+1)=0.
\end{equation*}
On the other hand suppose $\mathcal{H}^{s_{_0}}(E)>0$.  If for any
$\varepsilon>0$, $\mathcal{H}^{s_{_0}-\varepsilon}(E)$ were finite,
then by the above $\mathcal{H}^{s_{_0}}(E)=0$, a contradiction, whence
$\mathcal{H}^{s-\varepsilon}(E)=\infty$.
\begin{figure}[htb]
  \begin{center}
    \begin{picture}(300,220) 
      \put (40,10){\line(1,0){220}}
      \thicklines
      \put (40,200){\line(1,0){140}}
      \put (180,10){\vector(1,0){100}}
      \thinlines 
      \put (180,0){$s_0=\hdim{E}$}
      \put (40,0){$0$}
      \put (180,45){\circle*{3}}
      \put (270,0){$s$}
      \put (40,10){\line(0,1){115}}
      \multiput (40,120)(0,10){8}{\line(0,1){5}}
      \put (15,200){$\infty$}
      \put (40,200){\circle*{4}}
      \put (40,45){\circle*{3}}
      \multiput (180,10)(0,20){6}{\line(0,1){10}}
      \multiput (180,130)(0,10){7}{\line(0,1){5}}
      \put (0,120){$\mathcal H^s(E)$}
      \put (-5,45){$\mathcal H^{s_0}(E)$}
    \end{picture}
    \caption{ 
      The graph of the Hausdorff measure of $E$.  The Hausdorff
      dimension $\hdim{E}=s_0$, the point of discontinuity.  }
    \label{fig:HMgraph}
  \end{center}
\end{figure}
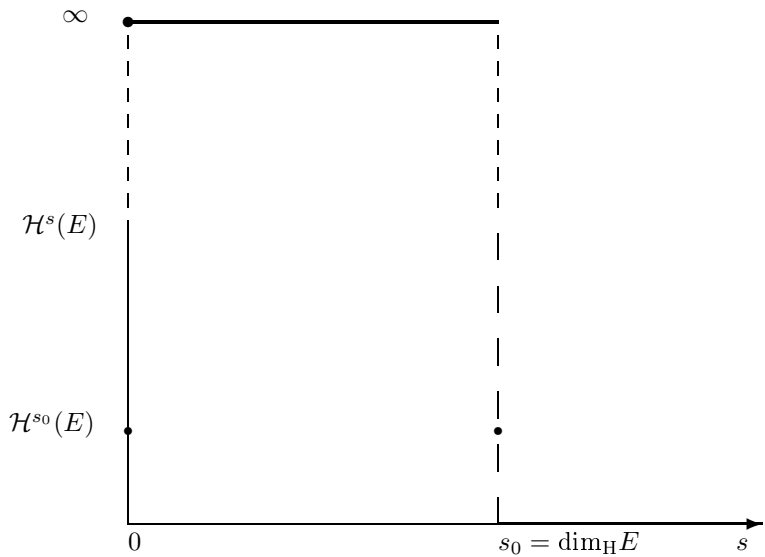

Thus for each infinite set $E$ in $n$-dimension\-al Euclidean space,
there exists a unique non-negative exponent $s_0$ such that
\begin{displaymath}
  \mathcal{H}^s(E)= 
  \begin{cases}
    \infty, 
    & \  0 \leqslant  s < s_0, \\
    0,  &    \  s_0 < s < \infty,   
  \end{cases}
\end{displaymath}
as shown in Figure~\ref{fig:HMgraph} (reproduced with the permission
of the Cambridge University Press from~\cite{MDAM}). The critical
exponent
\begin{equation}
  \label{eq:HDdef}
  s_0  = \inf\{s\in[0,\infty)\colon \mathcal{H}^s(E) = 0\}
\end{equation} 
where the Hausdorff $s$-measure crashes is called the \emph{Hausdorff
dimension} of the set $E$ and is denoted by $\hdim E$. Thus the
Hausdorff dimension of a finite set is 0, as it is for a countable
set.  It is clear that
\begin{equation}
  \label{eq:hdests}
  \text{if } \ \mathcal{H}^s(E)=0 \ \text{ then } \ \hdim E\leqslant
  s; \ \text{ and if } \mathcal{H}^s(E)>0 \ \text{ then } \ \hdim
  E\geqslant s.
\end{equation}

The Hausdorff dimension tells us nothing about the Hausdorff
$s$-measure at the critical exponent $s_0=\hdim E$, only that this is
the appropriate exponent to investigate the measure.  The sudden
change in Hausdorff $s$-measure at $s_0=\hdim E$ can be compared to the
focal length of a microscope.  If the lens is too close, the image
fills the eyepiece and cannot be resolved; if the lens is too far
away, the image is invisible.  At the focal length, the image is in
focus and can be resolved.

The main properties of Hausdorff dimension for sets in $\mathbb{R}^n$ 
are
\begin{enumerate}
\item[(i)] If $E \subseteq F$ then $\hdim E\leqslant \hdim F$.
\item[(ii)]  $\hdim E\leqslant n$.
\item[(iii)]  If  $\vert E\vert >0$, then $\hdim E = n$.\label{item:1}
\item[(iv)]  The dimension of a point is $0$.
\item[(v)] If $\hdim E < n$, then $\vert E\vert =0$ 
  (however $\hdim E =n$ does not imply $\vert E\vert >0$). 
\item[(vi)] $\hdim (E_1\times E_2) \geqslant \hdim E_1+\hdim E_2$
\item[(vii)] $\hdim \cup_{j=1}^\infty E_j=\sup\{\hdim E_j\colon 
  j\in\mathbb{N}\}$.
\end{enumerate}
It can be shown that the Hausdorff dimension of any countable set is 0
and that of any open set in $\mathbb{R}^n$ is
$n$~\cite[p.~29]{FalcFG}.  The nature of the construction of Hausdorff
measure ensures that the Hausdorff dimension of a set is unchanged by
an invertible transformation which is bi-Lipschitz.  This implies that
for any set $S \subseteq \mathbb{R} \setminus\{0\}$, $\hdim
S^{-1}=\hdim S$, where $S^{-1}=\{s^{-1}\colon s\in S\}$.  Thus on the
whole, Hausdorff dimension behaves as a dimension should, although the
natural formula 
\begin{displaymath}
  \hdim( E_1\times E_2)= \hdim E_1+\hdim E_2
\end{displaymath}
does \emph{not} always hold~\cite[\S5.3]{FalcGFSshort} (it does hold
for certain sets, \emph{e.g.}, cylinders, such as $E \times I$, where
I is an interval: $\hdim (E\times I)=\hdim E + \hdim I = \hdim E+1$ by
(iii), see~\cite{MDAM}).

The general character of $\delta$-covers in the definition of
Hausdorff outer measure can be difficult to work with and for many
applications in higher dimensions, it is convenient to restrict the
elements in the $\delta$-covers of a set to simpler sets such as balls
or cubes.  For example, covers consisting of hypercubes
\begin{displaymath}
  H=\{\mathbf{x}\in\mathbb{R}^n\colon
  \vert \mathbf{x}-\mathbf{a}\vert_\infty<\delta\},
\end{displaymath} 
where $\vert \mathbf{x}\vert_\infty:=\max\{\vert x_j\vert \colon
1\leqslant j\leqslant n\}$ is the {\em height} of
$\mathbf{x}\in\mathbb{R}^n$, centred at $\mathbf{a}\in\mathbb{R}^n$
and with sides of length $2\delta$ are used extensively.  While outer
measures corresponding to these more convenient restricted covers are
not the same as Hausdorff measure, they are comparable and so have the
same critical exponent~\cite[Chapter~5]{MattilaGS}.  Thus there is no
loss as far as dimension is concerned if the sets $U$ are chosen to be
balls or hypercubes. Of course, the two measures are identical for
sets with Hausdorff $s$-measure which is either $0$ or $\infty$. Such
sets are said to obey a `$0$-$\infty$' law, this being the appropriate
analogue of the more familiar `$0\,$-$1$' law in
probability~\cite[p~339]{KandT}.  Sets which do not satisfy a
$0$-$\infty$ law, \emph{i.e.}, sets which satisfy
\begin{equation}
  \label{eq:3}
  0<{\mathcal{H}}^{\hdim E}(E)<\infty,
\end{equation}
are called \emph{$s$-sets}; these occur surprisingly
often~\cite[p.~29]{FalcFG} and enjoy some nice properties (see
Chapters 2--4 of \cite{FalcGFSshort}).  One example is the Cantor set
which has Hausdorff $s$-measure 1 when $s=\log 2/\log
3$~\cite[p.~14]{FalcGFSshort}.  However it seems that $s$-sets are of
less interest in Diophantine approximation where the sets that arise
naturally, such as the set of badly approximable numbers or the set of
numbers approximable to a given order (see next section), obey a
$0$-$\infty$ law. The first
steps in this direction were taken by Jarník, who proved that the
Hausdorff $s$-measure of set of numbers rationally approximable to
order $v$ (see~\S\ref{sec:toa}) was $0$ or
$\infty$~\cite{Ja30a,Ja30b}.  This result turns on an idea related to
density of Hausdorff measure.
\begin{lemma}
  \label{lem:qis}
  Let $E$ be a null set in $\mathbb{R}$. Suppose that   for any interval 
  $(a,b)$
  and  $s\in [0,1]$,
  \begin{equation}
    \label{eq:qind}
    \mathcal{H}^s(E\cap(a,b)) \leqslant K(b-a) \mathcal{H}^s(E).  
  \end{equation}
  Then $\mathcal{H}^s(E)=0$ or $\infty$.
\end{lemma}
\begin{proof}
  Suppose the contrary, \emph{i.e.}, suppose $0 < \mathcal{H}^s(E)
  <\infty $.  Since $E$ is null, given $\varepsilon>0$, there exists a
  cover of $E$ by open intervals $(a_j,b_j)$ such that
  \begin{displaymath}
    \sum_j (b_j-a_j) < \varepsilon.
  \end{displaymath}
  By~\eqref{eq:qind}, there exists a constant $K>0$ such that
  \begin{equation*}
    0 < \mathcal{H}^s(E) = \mathcal{H}^s(\cup_j (a_j,b_j)\cap
    E) \leqslant K\mathcal{H}^s(E)\sum_j(b_j-a_j) <K \varepsilon
    \mathcal{H}^s(E)< \mathcal{H}^s(E)
  \end{equation*}
  for $\varepsilon<1/K$, a contradiction. 
\end{proof}

The proof for a general outer measure is essentially the same.  The
sets we encounter in Diophantine approximation are generally not
$s$-sets and some satisfy this `quasi-independence' property.

For other definitions of dimension, such as box-counting and packing
dimension, and their relationship with Hausdorff dimension,
see~\cite[Chapter~3]{FalcFG}.

\subsection{The determination of Hausdorff dimension}
Unless some general result is available, the Hausdorff dimension
$\hdim E$ of a null set $E$ is usually determined in two steps, with
the correct upward inequality $\hdim E\leqslant s_0$ and downward
inequality $\hdim E\geqslant s_0$ being established separately.

\subsubsection{The upper bound}
\label{sec:ub}
In view of~\eqref{eq:hdests}, an upper bound can be obtained by
finding a value of $s$ for which $\mathcal{H}^s(E)$ vanishes.  To find
such a value, it suffices to exhibit a cover $\{H\}$ of $E$
($E\subseteq \cup H$) by hypercubes $H$ of arbitrarily small
sidelength and $s$-length.  This can often be done by adapting the
estimate involved in showing that Lebesgue measure is 0. When $E$ is a
limsup set, \emph{i.e.},
\begin{equation*}
  E=\limsup_{N\to\infty}E_N=\bigcap_{N=1}^\infty 
  \bigcup_{k=N}^\infty
  E_k
  =\{\mathbf{x}\in\mathbb{R}^n\colon \mathbf{x}\in E_k  
  \mathrm{ \ for\ infinitely\ many \ } k\in \mathbb{N}\}
\end{equation*}
for a sequence of sets $E_n$, a simple Hausdorff measure counterpart
of the convergence case of the Borel-Cantelli lemma often gives the
correct upper bound for $\hdim E$.  This is useful in Diophantine
approximation.
\begin{lemma}
  \label{lem:HCLemma}  
  Let
  \begin{equation*}
    E= \{\mathbf{x}\in\mathbb{R}^n\colon \mathbf{x}\in E_k
    \mathrm{ \ for\ infinitely\ many \ } k\in \mathbb{N}\}. 
  \end{equation*}
  If for some $s>0$,  
  \begin{equation}
    \label{eq:slength}
    \sum_{k=1}^\infty \di(E_k)^s < \infty,   
  \end{equation}
  then $\mathcal{H}^s(E)=0$ and $\hdim E\leqslant s$.
\end{lemma}
\begin{proof}
  From the definition, for each
  $N=1,2,\dots$,
  \begin{equation*}
    E\subseteq \bigcup_{k=N}^\infty E_k,
  \end{equation*}
  so that the family $\scr{C}^{(N)}=\{E_k\colon k\geqslant N\}$ is a
  cover for $E$.  By~\eqref{eq:slength},
  \begin{displaymath}
    \lim_{N\to \infty} \sum_{k=N}^\infty \di(E_k)^s 
    =0.
  \end{displaymath}
  Hence $\lim_{k\to\infty} \di(E_k) = 0$ and therefore given
  $\delta>0$, $\scr{C}^{(N)}$ is a $\delta$-cover of $E$ for $N$
  sufficiently large. But 
  \begin{equation*}
    \mathcal{H}_\delta^s(E)
    =\inf_{\scr{C}_\delta} \sum_{U\in{\scr{C_\delta}}} (\di U)^{s}
    \leqslant   \ell^s(\scr{C}^{(N)})
    =\sum_{k=N}^\infty \di(E_k)^s  \to 0
  \end{equation*}
  as $N\to\infty$.  Thus $\mathcal{H}_\delta^s(E)=0$ and
  by~\eqref{eq:HHd}, $\mathcal{H}^s(E)=0$, whence $\hdim E\leqslant
  s$.
\end{proof}

\subsubsection{The lower bound}
\label{sec:lb}
The lower bound is often harder (though by no means always,
see~\cite{dd00}).  It requires showing that given \emph{any} $s<s_0$
and \emph{any} cover $\{C\}$ of $E$ with the diameters of the covering
elements arbitrarily small, the $s$-length $\sum_C (\di C)^s \geqslant
\delta$ for some positive $\delta$.  This can be very difficult and
has led to the development of a variety of methods.  In Diophantine
approximation, the \emph{regular systems} introduced by Baker and
W.~M.~Schmidt~\cite{BS} and the more general \emph{ubiquitous
systems}~\cite{DRV90a} depend on a good supply of approximating
elements (\emph{e.g.}, the rationals).  These and related techniques
have proved effective in obtaining lower bounds for the Hausdorff
dimension of sets of number theoretic interest (see the survey
article~\cite{bbdbaker} for more details). A more fundamental approach
is the so-called mass distribution principle.
\begin{lemma}
  Let $\mu$ be a measure supported on a bounded Borel set $E$ in
  $\mathbb{R}^n$. Suppose that for some $s\geqslant 0$, there are
  strictly positive constants $c$ and $\delta$ such that $\mu(B)
  \leqslant c\, (\di B)^s$ for any ball $B$ in $\mathbb{R}^n$ with
  $\di B\leqslant \delta$.  Then ${\mathcal{H}}^s(E)\geqslant
  \mu(E)/c$. 
\end{lemma}
\begin{proof}
  The proof is short.  Let $\{B_k\}$ be a $\delta$-cover of $E$ by
  balls $B_k$.  Then
  \begin{equation*}
    \mu(E)\leqslant \mu\left(\bigcup_k B_k\right) \leqslant \sum_k
    \mu(B_k) \leqslant c\sum_k (\di B)^s.
  \end{equation*}
  Taking infima over all such covers, we see that
  ${\mathcal{H}}_\delta^s(E)\geqslant \mu(E)/c$, whence on letting
  $\delta \to 0$,
  \begin{displaymath}
    {\mathcal{H}}^s(E)\geqslant \mu(E)/c>0.
  \end{displaymath}
\end{proof}
This simple lemma is surprisingly useful and gives the easy part of
Frostman's lemma~\cite{Frostman} which is now stated in full.  The
Vinogradov notation $a\ll b$ for $a,b>0$ means that $a=O(b)$.
\begin{lemma} 
  \label{lem:Frostman}
  Let $E$  be a Borel subset of $\mathbb{R}^n$.
  Then
  \begin{displaymath}
    \mathcal{H}^s(E)>0 
  \end{displaymath} 
  if and only if there exists a measure $\mu$ on $\mathbb{R}^n$
  supported on $E$ with $\mu(E)$ finite such that $\mu(B) \ll (\di B
  )^s$ for all sufficiently small balls $B$.
\end{lemma}
Thus if $E$ supports a probability measure $\mu$ ($\mu(E)=1$) with $
\mu(B)\ll (\di B )^s $ for all sufficiently small balls $B$, then
$\hdim E\geqslant s$. The converse is more difficult but can be proved
using net measures (see~\cite{CarlesonES,FalcFG,MattilaGS}).

\section{Diophantine approximation} 
At its simplest level, Diophantine approximation is concerned with
approximating real numbers by rationals. Hardy and Wright's classic
\emph{Introduction to the theory of numbers}~\cite{HW} contains an
excellent account while the more advanced~\cite{Casselshort,SchmidtDA}
are devoted wholly to Diophantine approximation.  The theory extends
to approximating vectors in $\mathbb{R}^n$ (simultaneous Diophantine
approximation) and to matrices (systems of linear forms).  For
simplicity, we will stick mainly to one particular direction in the
one dimensional real and complex cases and treat the extensions to
higher dimensions and other settings fairly briefly.  Since the
rationals $\mathbb{Q}$ are a dense subset of the real numbers
$\mathbb{R}$, given any real number $\alpha$ and any positive
$\varepsilon$, there exists a rational $p/q$ such that
\begin{equation}
  \label{eq:1}
  \left\vert \alpha - \dfrac{p}{q}\right\vert  < \varepsilon.
\end{equation}
The numerator $p$ is often of no interest and the size of the
expression
\begin{equation}
  \label{eq:4}
  \Vert q \alpha\Vert =\min\{\vert q \alpha -p\vert\colon 
  p\in\mathbb{Z}\},
\end{equation}
the distance of $q\alpha$ from the integers $\mathbb{Z}$, is
considered.  Although convenient, it will not be used much here in
order to keep the notational burden to a minimum.

In simultaneous Diophantine approximation, one considers the system of
$n$ inequalities
\begin{equation*}
  \left\vert \alpha_k-\dfrac{p_k}{q}\right\vert <\varepsilon, \ \ k=,1\dots,n.
\end{equation*}
This system can be expressed more concisely as a single vector
inequality with $\boldsymbol{\alpha}=(\alpha_1,\dots,\alpha_n)\in
\mathbb{R}^n$, $\mathbf{p}\in\mathbb{Z}^n$, $q\in\mathbb{N}$, by
considering the expression
\begin{equation*}
  \left\vert \boldsymbol{\alpha}-\dfrac{\mathbf{p}}{q}\right\vert_\infty ,
\end{equation*}
where for $\mathbf{x}\in\mathbb{R}^n$,
$\vert \mathbf{x}\vert_\infty=\max\{\vert x_1\vert ,\dots,\vert x_n\vert \}$ is the \emph{height}
of $\mathbf{x}$, or, on multiplying by $q$, the expression
\begin{equation*}
  \vert q\boldsymbol{\alpha}-\mathbf{p}\vert_\infty =\Vert
  q\boldsymbol{\alpha}\Vert =\max\{\Vert q\alpha_j\Vert\colon
  j=1,\dots,n\}.
\end{equation*}
The last inequality has a \emph{dual} or linear form version: given
$\boldsymbol{\alpha}\in\mathbb{R}^n$ and $\varepsilon>0$, one
considers the inequality
\begin{equation*}
  \vert\mathbf{q}\cdot\boldsymbol{\alpha}-p\,\vert<\varepsilon,
\end{equation*}
where $\mathbf{q}\in\mathbb{Z}^n$ and $p\in\mathbb{Z}$.
The last two inequalities can be combined into a single 
general one.  The
system of $n$ real linear forms
\begin{equation*}
  \xi_1 a_{1j} + \dots + \xi_m a_{mj}, \ j=1,\dots,n,
\end{equation*}
in $m$ real variables $\xi_1,\dots,\xi_n$, can be written more
concisely as $\boldsymbol{\xi} A$, where $A=(a_{ij})$ and the system
of $n$ inequalities in $m$ variables
\begin{align*}
  &\vert  q_1 a_{11} + \dots + q_m a_{m1} -p_1\vert < \varepsilon\\
  &\qquad \qquad\qquad \qquad \vdots \\
  & \vert  q_1 a_{1n} + \dots + q_m a_{mn} -p_n\vert < \varepsilon,
\end{align*}
can be written $\vert \mathbf{q} A-\mathbf{p}\vert_\infty <
\varepsilon$.  Further details are in~\cite{Casselshort,HW,SchmidtDA}.
The theory extends naturally to the fields of $p$-adic
numbers~\cite{CasselsLF,Lutz} and formal power
series~\cite{kristensen03,kristensen02}.  Less obviously, it also
extends to discrete groups acting on hyperbolic space.  This is
relevant to Diophantine approximation over the Gaussian integers or
rationals considered below in~\S\ref{sec:DAGI}, so an outline is now
given.

The hyperbolic space setting sprang from the observations that the the
real axis is the set of limit points of the rationals and that the
rationals can be characterised as the parabolic vertices of the
modular group $\SL(2,\mathbb{Z})$, \emph{i.e.}, as the orbit of the
point at infinity under the linear fractional or Möbius
transformations
\begin{equation}
\label{eq:Mobius}
  z\mapsto \dfrac{az+b}{cz+d}, \ a,b,c,d \in \mathbb{Z}, \ ad-bc=1,
\end{equation}
of the extended upper half plane
$\mathbb{H}^{\,2}=\{z=x+iy\in\mathbb{C}\colon
y\geqslant0\}\cup\{\infty\}$.  For each element
\begin{equation*}
  g= \left(\begin{matrix} a & b \\
      c & d\\
    \end{matrix}\right)
\end{equation*}
in $\SL(2,\mathbb{Z})$, the point $\infty$ in the extended real line
evidently goes to $a/c$ under the group action and is also a fixed
point of the map $z\mapsto z+1$.  The maps $g$ form the \emph{modular}
group, $\SL(2,\mathbb{Z})$, a discrete subgroup of
$\SL(2,\mathbb{R})$, which is essentially (modulo the centre) the
group of orientation preserving Möbius transformations of the upper
half plane $\mathbb{H}^{\,2}$ to itself.  When the upper half plane is
endowed with the hyperbolic metric derived from $d\rho = \abs{dz}/y$,
it is a model for two dimensional hyperbolic space
$(\mathbb{H}^{\,2},\rho)$. The Möbius group $M(\mathbb{H}^{\,2})$ is
the group of isometries of $(\mathbb{H}^{\,2},\rho)$. Because the
group $\SL(2,\mathbb{Z})$ is discrete, points in the orbit can
accumulate only on the boundary $\mathbb{R}\, \cup\, \{\infty\}$ of
$\mathbb{H}^{\,2}$ and because the group elements are isometries with
respect to the hyperbolic metric, the limit set of any orbit is
the extended real line $\mathbb{R}\, \cup\{\infty\}$. A discrete subgroup of
$M(\mathbb{H}^{\,2})$ is called a \emph{Fuchsian} group.  Further
details of this rich and beautiful theory are
in~\cite{AhlforsMT,AndersonHG,BeardonGDG,NichollsETDG,Patterson76a}
and there is a short survey in Chapter~7 of~\cite{MDAM}.

These observations allow the classical theory, including the metrical
theory, of Diophantine approximation to be translated into Fuchsian
groups acting on the hyperbolic plane and to the much more general
setting of \emph{Kleinian} groups acting on $(n+1)$-dimensional
hyperbolic space $(\mathbb{H}^{\,n+1},\rho)$, $n\geqslant 2$ (Kleinian
groups are the discrete subgroups of the M\"obius group
$M(\mathbb{H}^{\,n+1})$ of isometries of $\mathbb{H}^{\,n+1}$; further
details are in the references given). The \emph{Picard} group
$\SL(2,\mathbb{Z}[i])$ consists of $2\times 2$ matrices over
$\mathbb{Z}[i]$ with determinant 1 and has an action on
$\mathbb{C}$ given by 
\begin{equation}
\label{eq:Picard}
  z \mapsto  \frac{az + b}{cz+d},  
\ a,b,c,d \in \mathbb{Z}[i], \ ad-bc=1.
\end{equation}
The limit set of the Picard group is the extended complex plane (or
Riemann sphere) $\mathbb{C}\cup \{\infty\}$ and the orbit of the point
at $\infty$ under the group is the set of ratios of Gaussian integers
~\cite{MP93}.
Thus the Picard group plays a role precisely parallel to that played
by the modular group, expressed by~\eqref{eq:Mobius}, in approximating
real numbers by ratios of integers.

In the literature cited, hyperbolic space is usually taken in the
equivalent Poincar\'e form of the open unit ball
$\mathbb{B}^{n+1}=\{\mathbf{x}\in \mathbb{R}^{n+1} :
\abs{\mathbf{x}}_2<1\}$, where the ball is now endowed with the  
equivalent hyperbolic metric $\rho$ given by
$d\rho=\abs{d\mathbf{x}}_2/(1-\abs{\mathbf{x}}_2^2)$. To ease
comparison, we shall adopt this viewpoint, even though the upper half
plane model, as used by Sullivan in~\cite{Sullivan82}, is more
natural for Diophantine approximation. We choose to consider the
$(n+1)$-dimensional hyperbolic space, as the results in Diophantine
approximation are results about the boundary of $\mathbb{H}^{\,n+1}$,
which is $n$-dimensional.

The analogue $\mathfrak{p}$ of the point at infinity for Kleinian
groups is not quite straightforward. First of all the nature of the
elements $g$ of the Kleinian group $G$ implies that each $g$ has at
most two fixed points on the boundary of the ball. The special point
$\mathfrak{p}$ is called \emph{a parabolic} fixed point if it is the
unique fixed point on the boundary of some element in $G$; otherwise
they are called \emph{hyperbolic} fixed points. The orbit of a special
point $\mathfrak{p}$ under the action of a Kleinian group $G$
corresponds to the rationals $\mathbb{Q}$. The limit set $\Lambda(G)$
of the orbit under $G$ of a point in $\mathbb{H}^{\,n+1}$ lies in the
boundary $\mathbb{S}^n$. Given $\alpha\in
\Lambda(G)\subseteq\mathbb{S}^n$, one considers the 
quantity 
\begin{equation*}
  \vert \alpha- g(\mathfrak{p})\vert_2,
\end{equation*}
where $\vert \cdot \vert_2$ is the usual Euclidean metric in
$\mathbb{R}^{n+1}$.  Analogues of the principal theorems in Diophantine
approximation have been obtained with relatively minor technical
restrictions and will be discussed below; a brief survey is in
Chapter~7 of~\cite{MDAM}.  There is a striking dynamical
interpretation of the approximation in terms of flows on the
associated quotient space $\mathbb{H}^{n+1}/G$; more details are 
in~\cite{MDAM,FM95,MP93,Sullivan82,VelaniKTHD}.  
We now return to the one-dimensional theory.

\subsection{Dirichlet's theorem}

It is not difficult to make~\eqref{eq:1} more precise: given any real
number $\alpha$ and any positive integer $q$, there exists an integer $p$
such that $\vert q\alpha-p\vert  < 1$, and indeed such that
\begin{displaymath}
  \left\vert \alpha - \dfrac{p}{q}\right\vert  \leqslant \dfrac{1}{2q}.
\end{displaymath}
There are  denominators $q$ for which more can be said by
using Dirichlet's celebrated  `box argument' 
(see~\cite{Casselshort,HW}).
\begin{theorem}
  \label{thm:DT}
  For each real number $\alpha$ and any positive integer $N \geq 1$,
  there exists a rational $p/q$ with denominator satisfying
  $1\leqslant q \leqslant N$, such that
  \begin{equation*}
    \left\vert\alpha - \dfrac{p}{q}\right\vert < \dfrac{1}{qN}
    \leqslant \dfrac{1}{q^2}.
  \end{equation*}
\end{theorem}

\begin{proof}
  Let $[\alpha]$ be the integer part of $\alpha$ and $\{\alpha\}$ its
  fractional part, so that $\alpha=[\alpha]+\{\alpha\}$.  Divide the
  interval $[0,1)$ into $N$ subintervals $[k/N,(k+1)/N)$, where
  $k=0,1,\dots,N-1$, of length $1/N$.  The $N+1$ numbers $ \{r
  \alpha\}$, $r=0,1,\dots,N$, fall into the interval [0,1) and so two,
  $\{r \alpha\}$, $\{r'\alpha\}$ say, must fall into the same
  subinterval, $[k/N,(k+1)/N)$ say. Suppose that $r > r'$. Then
  \begin{displaymath}
    \vert \{r \alpha\} -\{r'\alpha\}\vert = \vert
    r\alpha-[r\alpha]-r'\alpha+[r'\alpha]\vert = \vert q\alpha-p\vert
    < \dfrac{1}{N},
  \end{displaymath}
  where $q=r-r', p = [r\alpha]-[r'\alpha] \in \mathbb{Z}$ and $1
  \leqslant q \leqslant N$.  Dividing by $q$ gives the quantitative 
  inequality
  \begin{equation}
    \label{eq:DTq}
    \left\vert\alpha - \dfrac{p}{q}\right\vert < \dfrac{1}{qN}
  \end{equation}
  and since $1\leqslant q\leqslant N$, the final inequality is
  immediate.
\end{proof}
A nice sharpening is in~\cite[p.~1]{HarmanMNT}.  When $p,q$ are
restricted to having highest common factor 1, the inequality
\begin{equation}
  \label{eq:DTio}
  \left\vert\alpha - \dfrac{p}{q}\right\vert < \dfrac{1}{q^2}
\end{equation}
has only finitely many solutions if and only if $\alpha$ is a
rational. Thus almost all real numbers satisfy~\eqref{eq:DTio} for
infinitely many rationals $p/q$.  Without this
restriction,~\eqref{eq:DTio} holds infinitely often for all
$\alpha\in\mathbb{R}$.

Dirichlet's theorem is one of the fundamental results in the theory of
Diophantine approximation. It can be viewed as a result about covers
and plays a central part in the Jarník--Besicovitch theorem, discussed
below. The theorem generalises to the simultaneous Diophantine
approximation of $n$ real numbers
$\alpha_1,\dots,\alpha_n$~\cite[Theorem~200]{HW} and asserts that 
given $N\in\mathbb{N}$, there exists $q\in\mathbb{N}$, with 
$q\leqslant N$ and $\mathbf{p}=(p_1,\dots,p_n)\in\mathbb{Z}^n$
such that
\begin{equation*}
  \left\vert \boldsymbol{\alpha}-\dfrac{\mathbf{p}}{q}\right\vert_\infty
  <\dfrac{1}{qN^{1/n}}.
\end{equation*}
In particular, for simultaneous Diophantine approximation in the
plane, given any $\boldsymbol{\alpha} = (\alpha_1, \alpha_2) \in
\mathbb{R}^2$ and $N \in \mathbb{N}$, there exists $q \in \mathbb{N}$ 
with $q \leqslant N$ and $p_1, p_2 \in \mathbb{Z}$ such that
\begin{equation}
  \label{eq:5}
  \max \left\{ \abs{\alpha_1 - \dfrac{p_1}{q}}, \abs{\alpha_2 -
  \dfrac{p_2}{q}} \right\} < \dfrac{1}{q N^{1/2}} \leqslant q^{-3/2}. 
\end{equation}

There is a so-called \emph{dual} version: given
$\boldsymbol{\alpha}\in\mathbb{R}^n$ and $N\in\mathbb{N}$, there
exists $\mathbf{q}\in\mathbb{Z}^n$ and a $p\in\mathbb{Z}$ such that
\begin{equation*}
  \vert \mathbf{q}\cdot\boldsymbol{\alpha}-p\vert <N^{-n}. 
\end{equation*}
These can be combined into a result for systems of $n$ real linear
forms~\cite[Chapter~1, Theorem~VI]{Casselshort}:
\begin{equation*}
  \vert \mathbf{q} A-\mathbf{p}\vert_\infty < N^{-m/n}.  
\end{equation*}

In the setting of a Kleinian group acting on hyperbolic space, the
analogue of the denominator in Theorem~\ref{thm:DT} corresponding to
$g(\mathfrak{p})$ is defined to be
\begin{equation}
  \label{eq:lg}
  \lambda_g:= \vert \det(Dg\vert_0)\vert ^{-1}=\dfrac{1}{2} \cosh 
  \rho(0,g(0))\asymp
  e^{\rho(0,g(0))}
\end{equation}
in the ball model, and so $\lambda_g\to\infty$ as $\vert g(0)\vert \to
1$, \emph{i.e.}, as the orbit of the origin moves towards the
boundary. Here $Dg\vert_0$ denotes the Jacobian of $g$ evaluated at
the origin. For finitely generated Fuchsian groups of the first kind
taken to be acting on the closed unit disc $\Delta$, the elements $g$
are of the form
\begin{equation*}
  g= \left(\begin{matrix} a & b \\
      \overline{b} & \overline{a}\\
    \end{matrix}\right), \ a,b \in \Delta, \vert a\vert^2- 
  \vert
  b\vert^2=1 
\end{equation*}
and $\lambda_g=2( \vert a\vert^2+ \vert b\vert^2)$. Hedlund's lemma,
which is a classical result in the theory of discontinuous groups, is a
partial analogue of Dirichlet's theorem.  Let $\zeta\in \mathbb{S}^1$,
the unit circle.  For any $\xi\in \mathbb{S}^1$ which is not a
parabolic point, there exist infinitely many $g\in G$ such that
\begin{equation*}
  \vert \xi-g(\zeta)\vert <\dfrac{C}{\lambda_g}\, ,
\end{equation*}
for some $C>0$.  A complete analogue of Dirichlet's theorem, including
the quantitative inequality~\eqref{eq:DTq}, was obtained by
Patterson~\cite{Patterson76a} for Fuchsian groups and later he and
others extended it to Kleinian
groups~\cite{Patterson76b,Stratmann94,SV95,VelaniKTHD}.  However, the
statements in the Kleinian group setting differ for parabolic and
hyperbolic fixed points and so for simplicity the result will be
stated when $G$ has a unique parabolic point $\mathfrak{p}$.  Let
$N\geqslant 2$. Then for any $\xi\in\Lambda(G)$, there exists a $g\in
G$ with $\lambda_g<N$ such that
\begin{equation*}
  \vert \xi- g(\mathfrak{p})\vert_2 <\dfrac{C}{\sqrt{N\lambda_g}},
\end{equation*}
where $C$ is a constant depending only on $G$. 

\subsection{Types of approximation}   
\label{sec:toa}
The equation~\eqref{eq:DTio} given by Dirichlet's theorem is
essentially best possible as Hurwitz~\cite{Hurwitz1891} showed that
each $\alpha \in \mathbb{R}$ satisfied the inequality
\begin{equation*}
  \left\vert \alpha-\dfrac{p}{q}\right\vert  < \dfrac{1}{\sqrt{5}q^2}
\end{equation*}
for infinitely many positive integers $q$ and that this was best
possible in the sense that the constant $1/\sqrt{5}$ cannot be reduced
for numbers $\alpha$ equivalent to the golden ratio $(\sqrt{5}\,+1)/2$
\cite{Casselshort,Forder63,HW,SchmidtDA}.

\subsubsection{Badly approximable numbers} 
\label{sec:BA}
A number $\beta\in \mathbb{R}$ is called \emph{badly approximable} or
of \emph{constant type} if there exists a $K=K(\beta)$ such that
\begin{equation}
  \label{eq:ba}
  \left\vert \beta-\dfrac{p}{q}\right\vert \geqslant \dfrac{K}{q^2} 
\end{equation}
for all $p/q\in\mathbb{Q}$. Using \eqref{eq:4}, \eqref{eq:ba} can be
written as $q \Vert q\beta \Vert \geqslant K$. In view of Hurwitz's
theorem, the constant $K<1/\sqrt5$. Quadratic irrationals, such as
$\sqrt2$ and the golden ratio $(\sqrt{5}+1)/2$, are badly
approximable. This is proved in~\cite[\S11.4]{HW} using the fact that
the partial quotients in the continued fraction expansion for a
quadratic irrational are periodic. However, the proof relies only on
the boundedness of the partial quotients, which therefore
characterises the badly approximable numbers. The set of badly
approximable numbers will be denoted by $\mathfrak{B}$. The notion
extends naturally to higher dimensions and to the more general
settings mentioned above. 

Badly approximable numbers are important in applications, particularly
in stability questions for certain dynamical systems
\cite{MR93g:11011}.  For example the `noble' numbers, which are
equivalent to the golden ratio, have been conjectured to be the most
robust in the breaking up of invariant tori~\cite{Mackay83}. One very
practical application involved the design of rocket casings. These
were made using ruled surfaces and vibrations from the motors were
propagated along the generators. To reduce the effects of resonance
and delay the onset of catastrophic vibration, the ratio of the
circumference to the length of the casing was chosen to be a quadratic
irrational (V.~I.~ Arno'ld, personal communication). The desirable
properties of badly approximable numbers (and in particular of the
golden ratio) appear to be related to their occurrence in nature. It
has recently been discovered that the ratio between two step heights
on the surface of certain quasi-crystals is given by the golden ratio
(see~\cite{chown02:_why} for statements of this result and additional
examples).

The notion of badly approximable numbers carries over to higher
dimensions, including systems of linear forms~\cite{Schmidt69},
$p$-adics~\cite{AA95a} and fields of formal power series
\cite{kristensen03} as well as to Kleinian groups acting on hyperbolic
space~\cite{BJ97,FM95,Patterson76a,Patterson89}.  In the hyperbolic
space setting, a point $\beta$ in $\Lambda(G)$ is said to be badly
approximable with respect to $\mathfrak{p}$ if there exists a positive
constant $K=K(\beta)$ such that
\begin{equation*}
  \abs{\beta - g(\mathfrak{p})}_2\geqslant K/\lambda_g
\end{equation*}
for all $g\in G$.

\subsubsection{Diophantine type}
\label{sec:Dt}
The concept of a badly approximable number has extensions to
restricted classes of real numbers and points in $\mathbb{R}^n$ that are 
useful in connection with stability and other questions
(see~\S\ref{sec:applns}) and fortunately enjoys full measure. Let
$K>0$, $v>1$.  The real number $\alpha$ is said to be of
\emph{Diophantine type} $(K,v)$~\cite{ArnoldMM} (the definition has
been altered slightly for consistency) if
\begin{equation*}
  \left\vert \alpha-\dfrac{p}{q}\right\vert \geqslant \dfrac{K}{q^{v+1}} 
\end{equation*}
for all rationals $p/q$; the set of numbers of Diophantine type
$(K,v)$ is denoted by $\mathcal{D}(K,v)$. The union
\begin{equation*}
  \mathcal{D}_v=\bigcup_{K>0}\mathcal{D}(K,v)
\end{equation*}
consists of numbers of Diophantine type $v$ (\emph{i.e.}, of type
$(K,v)$ for some $K>0$) and the union
\begin{equation*}
  \mathcal{D}= \bigcup_{v>1} \mathcal{D}_v= \bigcup_{K>0,\, v>1}\mathcal{D}(K,v)
\end{equation*}
is the set of numbers of Diophantine type $v$ for some $v>1$.  Note
that $\mathcal{D}_1=\mathfrak{B}$, the set of badly approximable
numbers.  We will see in~\S\ref{sec:ktamda} that $\mathfrak{B}$ is
null but that when $v>1$, $\mathcal{D}_v$ has full measure, which is
pleasing since points of Diophantine type have desirable approximation
properties for certain applications.  Again, the notion of Diophantine
type extends naturally to higher dimensions~\cite{DRV94}.  We will be
particularly interested in the case of a single linear form and
accordingly we extend the definition of Diophantine type for a real
number to a point in $\mathbb{R}^n$. A point
$\boldsymbol{\beta}\in\mathbb{R}^n$ is of (dual) Diophantine type
$(K,v)$ if for all $p\in\mathbb{Z}$ and non-zero
$\mathbf{q}\in\mathbb{Z}^n$,
\begin{equation*}
  \vert \mathbf{q}\cdot \boldsymbol{\beta}-p\vert \geqslant
   \dfrac{K}{\vert \mathbf{q}\vert_\infty^{v}}.
\end{equation*}

\subsubsection{Well approximable numbers} 
\label{sec:wan} 

In applications, we will be interested in points $\alpha$ which are
not of Diophantine type $(K,v)$ for any $K>0$, \emph{i.e.}, in one
dimension with the set
\begin{equation}
  \label{eq:vwa1}
  E_v=\left\{\alpha\in\mathbb{R}\colon  {\text{ for any }} \ K>0, 
    \left\vert  \alpha-\dfrac{p}{q}\right\vert  <
    \dfrac{K}{q^{v+1}}  {\text{ for some }} \ \dfrac{p}{q}\in\mathbb{Q}
  \right\}=\mathbb{R}\setminus \mathcal{D}_v,
\end{equation}
the complement of $\mathcal{D}_v$.  These numbers are closely related
to numbers which are \emph{rationally approximable to order
  $v+1$}~\cite[\S11.4]{HW}, \emph{i.e.}, to the set
\begin{equation}
  \label{eq:Rv}
  R_v=\left\{\alpha\in\mathbb{R}\colon  {\text{ for some }} \ K>0, 
    \left\vert  \alpha-\dfrac{p}{q}\right\vert  <
    \dfrac{K}{q^{v+1}}  {\text{ for infinitely many }} \ 
    \dfrac{p}{q}\in\mathbb{Q} \right\}
\end{equation}
(the exponent $v+1$ in the denominator is a normalisation to keep the
notation the same as elsewhere.) By Dirichlet's theorem, all real
numbers are rationally approximable to order $2$ and quadratic
irrationals are rationally approximable to order exactly $2$.

The constant $K$ in the definitions of $R_v$ and $E_v$ is of less
significance than the exponent, which motivates the next definition.
A number $\alpha$ which satisfies
\begin{equation}
  \label{eq:Wvineq}
  \left\vert  \alpha-\dfrac{p}{q}\right\vert  < \dfrac{1}{q^{v+1}}
\end{equation}
for infinitely many $ p/q\in \mathbb{Q}$ will be called {\em
  $v$-approximable}; if $v>1$, $\alpha$ is called \emph{very well
  approximable}.  Thus Liouville numbers, which
satisfy~\eqref{eq:Wvineq} for any $v$, are very well approximable.
We let $W_v$ denote the limsup set of $v$-approximable numbers,
\emph{i.e.},
\begin{equation*}
  \label{eq:Wv}
  W_v=  \left\{\alpha\in\mathbb{R}\colon   
    \left\vert  \alpha-\dfrac{p}{q}\right\vert  <
    \dfrac{1}{q^{v+1}}  {\text{ for infinitely many }} \
    \dfrac{p}{q}\in\mathbb{Q} \right\}. 
\end{equation*}

The sets $E_v$, $R_v$, $W_v$ are related and decrease as $v$
increases, as follows from the inclusions
\begin{equation}
  \label{eq:UEWincl}
  R_{v+\varepsilon} \subset E_v\subset W_v\subset   R_{v},
\end{equation}
where $\varepsilon>0$ is arbitrary.  Clearly $W_v\subset R_v$.  Next
consider the complementary inclusion $\mathbb{R}\setminus W_v\subset
\mathcal{D}_v$ and let $ \alpha \not \in W_v$.  Then for all but
finitely many positive integers $q$,
\begin{equation*}
  \left\vert \alpha-\dfrac{p}{q}\right\vert \geqslant
  \dfrac{1}{q^{v+1}}, \quad j=1,\dots, r.
\end{equation*}
Moreover, for the exceptional values of $q$, say
$q^{(1)},\dots,q^{(r)}$, 
\begin{equation*}
  K=\min\left\{\abs{q^{(j)}\alpha-p^{(j)}}\left(q^{(j)}\right)^v\colon 
    j=1,\dots,r\right\} >0,
\end{equation*}
since $ q\alpha\in \mathbb{Z}$ implies $ k q\alpha\in\mathbb{Z}$ for
each $k\in\mathbb{Z}$.  Thus $\vert \alpha-p/q\vert  \geqslant
K/q^{v+1}$ for all $p/q\in \mathbb{Q}$ and so
$\alpha\in\mathcal{D}_v$, \emph{i.e.}, $\alpha\notin E_v$.

To establish the final inclusion, let $\alpha\in R_{v+\varepsilon}$,
so that for some $K>0$,
\begin{equation*}
  \left\vert\alpha-\dfrac{p}{q}\right\vert <
  \dfrac{K}{q^{v+1+\varepsilon}}
\end{equation*}
for infinitely many $p/q\in\mathbb{Q}$. Given any $K'>0$, choose as we 
may, a denominator $q_0$ sufficiently large so that
$q_0^{-\varepsilon}\leqslant K'/K$. Then there exists a
$p_0\in\mathbb{Z}$ such that
\begin{equation*}
  \left\vert \alpha-\dfrac{p_0}{q_0}\right\vert< 
  \dfrac{K'}{q_0^{v+1}}
\end{equation*}
and $\alpha\in E_v$.  The Jarník--Besicovitch
theorem~(\S\ref{sec:JBt}) gives us the Hausdorff dimension of $W_v$
and allows us to deduce from~\eqref{eq:UEWincl} that all the sets have
the same Hausdorff dimension.

\subsubsection{$\Psi$-approximable  numbers} 
\label{sec:psiapp}
Now we look at the set of points which enjoy a more general
approximation by rationals.  A function $\Psi\colon
\mathbb{N}\to\mathbb{R}^+$ such that
\begin{equation*}
  \label{eq:appfn}
  \lim_{q\to\infty}\Psi(q)=0 
\end{equation*}
will be called an \emph{approximation function}; without loss of
generality we can take $\Psi(q)\leqslant 1/(2q)$ and later we shall
also assume that $q\Psi(q)$ is decreasing (by decreasing we mean
non-increasing here and subsequently). A real 
number $\alpha$ is said to be \emph{$\Psi$-approximable} if $\alpha$
satisfies the inequality
\begin{equation}
  \label{eq:1d}  
  \left\vert \alpha-\dfrac{p}{q} \right\vert  < \Psi(q)
\end{equation}
for infinitely many $p/q\in\mathbb{Q}$ (there should be no confusion
with $v$-approximable numbers defined above). Note that this should
not be confused with the inequality $\vert q\alpha-p\vert <\psi(q)$
which is often considered, particularly in higher dimensions and that
there are other definitions depending on the form of~\eqref{eq:1d},
see for example \cite{bkm01}. The set of $\Psi$-approximable numbers
in $\mathbb{R}$ will be denoted $W(\Psi)$.

Since for each $k\in \mathbb{Z}$, $(p+kq)/q \in \mathbb{Q}$ and 
\begin{equation*}
  \left\vert \alpha+k-\dfrac{p+kq}{q} \right\vert =\left\vert
  \alpha-\dfrac{p}{q} \right\vert ,
\end{equation*}
it follows that $ W(\Psi)\cap[k,k+1)=(W(\Psi)\cap[0,1))+k$, so that
$W(\Psi)$ can be decomposed into a union over unit intervals:
\begin{displaymath}
  W(\Psi)=\bigcup_{k\in\mathbb{Z}} (W(\Psi)\cap[k,k+1))
  =\bigcup_{k\in\mathbb{Z}}  (W(\Psi)\cap[0,1))+k .
\end{displaymath}

As with the other types of approximation, the definitions of very well
approximable, $v$-approximable and $\Psi$-approximable numbers extend
naturally to systems of linear forms, $p$-adic numbers, formal power
series and to the hyperbolic setting. We will be interested in the
case of a single linear form: we say with some abuse of notation that 
a vector $\boldsymbol{\alpha}\in\mathbb{R}^n$ is $v$-approximable if
\begin{equation}
  \label{eq:valf}
  \vert \mathbf{q}\cdot \boldsymbol{\alpha}-p\vert 
  <\dfrac{1}{\vert \mathbf{q}\vert_\infty^{v}} 
\end{equation}
holds for infinitely many $\mathbf{q}\in\mathbb{Z}^n$ and
$p\in\mathbb{Z}$; $\boldsymbol{\alpha}$ is very well approximable if
$v>n$.  In the hyperbolic space setting, a point $\boldsymbol{\alpha}$
in $\mathbb{H}^{\,n}$ is $\Psi$-approximable with respect to $G$ and
$\mathfrak{p}$ if
\begin{displaymath}
  \left\vert\boldsymbol{\alpha}-g(\mathfrak{p})\right\vert_2
  <\Psi(\lambda_g)
\end{displaymath}
for infinitely many $g\in G$. These and the other analogues of the
real one-dimensional case will be discussed further when they arise
below.

\section{Khintchine's  theorem and metrical Diophantine approximation}
\label{sec:ktamda}
Following earlier work of Borel~\cite{Borel12}, Khintchine gave an
almost complete answer to the solubility of~\eqref{eq:1d} in terms of
the measure of $W(\Psi)$.  In a series of papers in the 1920's on the
Lebesgue measure of the sets $W(\Psi)$ and
$\mathfrak{B}$~\cite{Kh24,Kh25,Kh25a,Kh26}, he laid the foundations of
metrical Diophantine approximation.  This theory, which is closely
related to probability, measure and ergodic theory, considers sets of
solutions to Diophantine inequalities in terms of Lebesgue and other
measures. As a result, $0$-$1$ laws are a feature of the Lebesgue part
of the theory, as in Khintchine's theorem below (see
also~\cite[\S2.2]{HarmanMNT}). In addition, because an exceptional set
for which a result is invalid can be of measure zero, this can lead to
theorems having a strikingly simple yet general character.
\begin{theorem}[Khintchine]
  \label{thm:K1}
  \begin{equation*}
    \vert W(\Psi)\cap[0,1)\vert = 
    \begin{cases} 
      0 , &\text{if }  \sum_{k=1}^\infty k\Psi(k)<\infty, \\
      1, &\text{if } \ k\Psi(k) \ \text {is decreasing and }
      \sum_{k=1}^\infty k\Psi(k)=\infty. 
    \end{cases}
  \end{equation*}
\end{theorem}
Thus $W(\Psi)$ is null when the series $\sum_k k\Psi(k)$ converges and
is full when $k\Psi(k)$ decreases and the sum diverges. 
Subsequently Khintchine extended the result to simultaneous
Diophantine approximation~\cite{Kh26} (see also~\cite{Casselshort})
and Groshev extended it to systems of linear forms~\cite{Sprindzuk}.
In particular the measure of the set of points $(\alpha_1, \alpha_2)
\in [0,1)^2$ such that 
\begin{displaymath}
  \max\left\{ \abs{\alpha_1 - \dfrac{p_1}{q}}, \abs{\alpha_2 -
      \dfrac{p_2}{q}} \right\} < \Psi(q) 
\end{displaymath}
for infinitely many $p_1, p_2\in \mathbb{Z}$, $q \in \mathbb{N}$ is
$0$ or $1$ accordingly as $\sum k^2 \Psi(k)^2$ converges or as $k
\Psi(k)$ decreases and the sum diverges (\emph{cf.} Theorem
\ref{thm:complexKT} below).

The convergence case follows readily from the Borel-Cantelli
lemma~\cite{KandT} and the argument is closely related to that of
Lemma~\ref{lem:HCLemma}. However, the case of divergence is much more
difficult and relies on a crucial `pairwise quasi-independence' result
combined with `mean and variance' ideas or with results from ergodic
theory.  For further details of this remarkable improvement of the
Borel-Cantelli lemma,
see~\cite{Chung,ChungErdos52,HarmanMNT,KS64,RenyiPT,Sprindzuk,Sullivan84}
(the books~\cite{HarmanMNT,Sprindzuk} also include accounts of
W.~M.~Schmidt's important quantitative extension of Khintchine's
theorem~\cite{Schmidt60}).  In terms of probability, we recall that
the events $E_j$, $j=1,\dots,\infty$, with corresponding probabilities
$P(E_j)$ are pairwise independent if for any $j\ne k$, the probability
of the two events occurring is given by
\begin{displaymath}
   P(E_j\cap E_k)=P(E_j)P(E_k).
\end{displaymath}
As a result, the divergence half of the Borel-Cantelli Lemma holds for
pairwise independence.  Total independence, where the probability of
any finite sequence of events is given by the product of the
individual probabilities, is not necessary. Recently it has been shown
that the divergence case is also related to the lower bound for
Hausdorff dimension (see \S\ref{sec:cJBt}).

There has been dramatic progress in the metrical theory over the last
decades, particularly in the theory of `dependent variables' where the
point $\boldsymbol{\alpha}$ lies on a manifold, so that coordinates
of the point are functionally related. Sprind{\v z}uk's
solution~\cite{Spr80} of Mahler's conjecture in transcendence theory
in terms of the Diophantine approximation of points on the Veronese
curve $\{(x,x^2,\dots,x^n)\colon x\in\mathbb{R}\}$ gave this topic an
enormous impetus which has seen the recent proof of Sprind{\v z}uk's
conjectures and further results~\cite{bbdd99,bkm01,KM98a}. 

The $p$-adic analogue of Khintchine's theorem was obtained by
Jarník~\cite{Ja45} and extended to systems of linear forms by
Lutz~\cite{Lutz}. For fields of formal power series, the analogue was
obtained by de~Mathan~\cite{deMa70} and has recently been extended to
systems of linear forms~\cite{kristensen02}.  The complex analogue of
Khintchine's theorem, in which the approximation is by ratios of
Gaussian integers, is discussed in~\S\ref{sec:DAGI} below.

Khintchine's theorem corresponds to our intuition since if the
approximation function $\Psi$ is large, then there is a better chance
of the inequality being satisfied. In particular, the set $W_v$ is
null for $v>1$, since the series $\sum_k k^{-v}$ converges, and full for
$v\leqslant 1$ (when $W_v=[0,1)$ by Dirichlet's theorem).
\begin{cor}
  \begin{equation*}
    \vert W_v\cap[0,1)\vert = \vert R_v\cap[0,1)\vert =\vert
    E_v\cap[0,1)\vert = 
    \begin{cases}
      0, & \text{when } \ v>1,\\
      1, & \text{when }\ v\leqslant 1.
    \end{cases}
  \end{equation*}
\end{cor}
The result for the other sets follows from~\eqref{eq:UEWincl}.  Less
obviously, the theorem shows the Lebesgue measure of the set of
$\alpha\in [0,1)$ such that~\eqref{eq:1d} has infinitely many
solutions is $1$ when $\Psi(k) = 1/(k^2 \log k)$ and $0$ when $\Psi(k)
= 1/(k^2(\log k)^{1+\varepsilon})$ for any positive $\varepsilon$.

As a consequence of Jarník's theorem on simultaneous Diophantine
approximation to be discussed in~\S\ref{sec:JBt}, none of the sets
$W_v,R_v,E_v$ is an $s$-set (see~\eqref{eq:3} for the definition).
However, using Lemma~\ref{lem:qis} and invariance under rational
translates of $R_v$, Jarník had shown earlier that the set $R_v$ of
numbers rationally approximable to order $v$ is not an $s$-set and
obeys a `$0$-$\infty$' law~\cite{Ja30a,Ja30b}. Although superseded by
the above theorem, the argument is very nice but the papers cited are
not readily available, so the proof is repeated here.  Let $k$ be a
positive integer and let $\alpha\in R_v\cap[0,1/k]$, so that for some
$K=K(\alpha)>0$,
\begin{displaymath}
  \left\vert \alpha - \dfrac{p}{q}\right\vert<\dfrac{K}{q^{v+1}}
\end{displaymath}
for infinitely many rationals $p/q$.  For each $j/k$, $1\leqslant
j \leqslant k-1$, 
\begin{equation*}
  \left\vert\alpha - \dfrac{p}{q}\right\vert=  
  \left\vert \alpha+\dfrac{j}{k} - \dfrac{p}{q}-\dfrac{j}{k}\right\vert 
  < \dfrac{K k^{v+1}}{(kq)^{v+1}}.
\end{equation*}
Thus $\alpha+j/k\in R_v\cap[j/k,(j+1)/k]$ and it follows that for each
$j=0,1,\dots,k-1$,
\begin{displaymath}
  R_v\cap \left[\dfrac{j}{k},\dfrac{j+1}{k}\right] = 
  R_v\cap \left[0,\dfrac{1}{k} \right]+\dfrac{j}{k}.  
\end{displaymath}
Moreover since Hausdorff measure is translation invariant,
\begin{displaymath}
  \mathcal{H}^s\left(R_v\cap [0,1]\right)
  =\sum_{j=0}^{k-1}\mathcal{H}^s\left(R_v\cap\left[
  \dfrac{j}{k},\dfrac{j+1}{k} \right]\right)=
  k \mathcal{H}^s\left(R_v\cap\left[0,\dfrac{1}{k}\right]\right),
\end{displaymath}
whence for any $j,k$ with $1\leqslant j \leqslant k-1$,
\begin{displaymath}
  \mathcal{H}^s\left(R_v\cap\left[\dfrac{j}{k},\dfrac{j+1}{k}
    \right]\right)
  =\dfrac{1}{k}\mathcal{H}^s\left(R_v\cap[0,1])\right).
\end{displaymath}
Now every open interval $(a,b)$ can be represented as a union of a
countable set of intervals $[j/k,(j+1)/k]$, \emph{i.e.},
\begin{displaymath}
  (a,b) = \bigcup_{j,k} \left[\dfrac{j}{k},\dfrac{j+1}{k}\right],  
\end{displaymath}
so that $b-a = \sum_{j,k} 1/k$.  Hence
\begin{displaymath}
  \mathcal{H}^s((a,b)) = \mathcal{H}^s\left(\bigcup_{j,k} 
  \left[\dfrac{j}{k},\dfrac{j+1}{k}\right]\right)  
\end{displaymath}
and since $\mathcal{H}^s(\cdot)$ is an outer measure,
\begin{align*}
  \mathcal{H}^s\left((a,b)\cap R_v\cap[0,1]\right) &= 
  \mathcal{H}^s\left(\bigcup_{j,k} 
  R_v\cap \left[\dfrac{j}{k},\dfrac{j+1}{k}\right]\right)
  \leqslant \sum_{j,k} 
  \mathcal{H}^s\left(R_v\cap\left[\dfrac{j}{k},\dfrac{j+1}{k}
    \right]\right) \\ 
&  \leqslant \mathcal{H}^s(R_v\cap[0,1])\sum_{j,k} \dfrac{1}{k}
  =  (b-a)\mathcal{H}^s(R_v\cap[0,1]).
\end{align*}
Thus the hypotheses of Lemma~\ref{lem:qis} are satisfied and so
$\mathcal{H}^s(R_v\cap[0,1])=\mathcal{H}^s(R_v)$ is $0$ or
$\infty$. The argument can be extended to show that the set of numbers 
approximable to order $v$ by algebraic irrationals is not an
$s$-set. The measure at the critical exponent was shown to be $\infty$
by Bugeaud \cite{MR1931199}.

Khintchine's theorem also implies that the set $\mathfrak{B}$ of badly
approximable numbers is null. For given any $K>0$, the sum $\sum_q
(K/q)$ diverges and so by Khintchine's theorem the set of real numbers
$\alpha$ satisfying $ \vert \alpha-p/q\vert  < K/q^2 $ for infinitely
many $p/q\in \mathbb{Q}$ is full.  Thus the complementary set $F(K)$
of the set of $\alpha$ such that $\vert \alpha-p/q\vert  \geqslant
K/q^2$ for all but finitely many $p/q$ is null and evidently increases
as $K$ decreases. From its definition,
\begin{equation*}
  \mathfrak{B}\subset \bigcup_{K>0} F(K) = \bigcup_{N=1}^\infty
  F(1/N),
\end{equation*}
a countable union of null sets, whence $\mathfrak{B}$ is null. 

Since the set $W_v$ of very well approximable numbers is null, the
inclusions~\eqref{eq:UEWincl} imply that $\mathcal{D}_v$ is full for
$v>1$, \emph{i.e.}, that almost all real numbers are of Diophantine
type $(K,v)$ for some positive $K$.  Thus almost all numbers are
neither well or badly approximable. It is a remarkable fact, proved by
Jarník~\cite{Ja28} in 1928, the same year as Besicovitch's first paper
on Hausdorff measure and dimension (on planar 1-sets)~\cite{Bes28},
that although the set $\mathfrak{B}$ is null, $\hdim \mathfrak{B}=1$,
\emph{i.e.}, its Hausdorff dimension is maximal in the sense that it
coincides with that of the ambient space $\mathbb{R}$.  This we now
discuss.

\subsection{Jarn\'\i k's theorem for badly approximable 
  numbers}
\label{sec:Jathmba}

Let $\theta\in (0,1)$ and let $a_n$, $n=1,2,\dots$, denote the partial
quotients of the continued fraction expansion for $\theta$.  For each
$N\in\mathbb{N}$, define
\begin{equation*}
  M_N:=\{\theta\colon a_n\leqslant N\}.
\end{equation*}
Now a number is badly approximable if and only if it has bounded
partial quotients $a_n$~\cite[\S11.4]{HW}, so that
$\mathfrak{B}=\lim_{N\to\infty} M_N$.  In a pioneering paper that was 
the first on Hausdorff dimension in Diophantine approximation,
Jarník~\cite{Ja28} showed that for each $N\geqslant 8$,
\begin{equation*}
  1-\dfrac{4}{N\log 2}\leqslant \hdim M_N \leqslant 1-\dfrac{1}{8N\log
    N}.
\end{equation*}
It is of course immediate that $\hdim \mathfrak{B}\leqslant 1$ since
$\mathfrak{B}\subset \mathbb{R}$. We now state Jarník's theorem for
badly approximable numbers, which follows from the above.
\begin{theorem}[Jarník]
  \begin{displaymath}
    \hdim \mathfrak{B} = 1.
  \end{displaymath}
\end{theorem}

Using $(\alpha,\beta)$ games, W.~M.~Schmidt~\cite{Schmidt66,Schmidt69}
proved much more, extending Jarník's theorem to higher dimensions, so
that the set of simultaneously badly approximable points in the plane
has Hausdorff dimension $2$. In fact, $\mathfrak{B}$ is a \emph{thick}
set.  This is a `local' property in the sense that for each open
interval $I$, $\mathfrak{B}$ has `full' Hausdorff dimension,
\emph{i.e.}, $\hdim \mathfrak{B} \cap I=1$. A very general
inhomogeneous analogue has been proved using quite different ideas
from dynamical systems~\cite{Kleinbock99}. 

The game in one dimension involves two players A and B, a non-empty
set $S\subset \mathbb{R}$ and two parameters $\alpha\in (0,1)$, given
to the player A, and $ \beta\in (0,1)$, given to the player B.  Player
B begins by picking a closed interval $B_1$. Then A chooses a closed
subinterval $A_1\subset B_1$ with $\vert A_1\vert =\alpha\vert
B_1\vert $.  Then B picks an interval $B_2\subset A_1$ with $\vert
B_2\vert =\beta\vert A_1\vert =\beta\alpha \vert B_1\vert $ 
and then A chooses another subinterval $A_2\subset B_2$ with
$\vert A_2\vert =\alpha \vert B_2\vert =\alpha^2\beta\vert B_1\vert $
and so on.  Clearly the intervals $B_1, A_1, B_2, A_2, \dots,$ form a
decreasing nested sequence so that their intersection is a point,
$\omega$ say, in $B_1$.  Player A is called the winner if $\cap_j 
A_j=\{\omega\}\subseteq S$, otherwise B wins.

A fuller account of the game in $\mathbb{C}$ will be given in
\S\ref{sec:cba}, so we will simply say that Schmidt showed that when
$S=\mathfrak{B}$, A can force $\omega$ to be badly approximable, even
though $\mathfrak{B}$ is null, and deduced that $\hdim
\mathfrak{B}\geqslant 1$.  Since $\mathfrak{B}\subset \mathbb{R}$,
$\hdim \mathfrak{B}=1$.  Note that since
$\mathcal{H}^1(\mathfrak{B})=\vert \mathfrak{B}\vert $, the Lebesgue
measure of $\mathfrak{B}$ and since $\vert \mathfrak{B}\vert =0$
(Khintchine's theorem), it follows that the Hausdorff measure of
$\mathfrak{B}$ vanishes at $s=\hdim \mathfrak{B}$ and
\begin{equation*}
  \mathcal{H}^s(\mathfrak{B}) = 
  \begin{cases}
    \infty, & \   0 \leqslant  s < 1,\\ 
    0,  &    \  s\geqslant1.
  \end{cases}
\end{equation*} 
Thus $\mathfrak{B}$ is not an $s$-set (this can also be proved using
Jarník's lemma in \S\ref{sec:HD} above).

Jarník's theorem has also been extended to $p$-adic
fields~\cite{AA95a}, to fields of formal power 
series~\cite{kristensen03} and to hyperbolic space.  It follows from
the hyperbolic space counterpart of Khintchine's theorem that for
geometrically finite Kleinian groups $G$, the set
$\mathfrak{B}(G,\mathfrak{p})$ of the hyperbolic analogue of badly
approximable points, has zero Patterson measure and the analogue of
Jarník's theorem holds~\cite{FM95,Patterson76a,Patterson89},
\emph{i.e.},
\begin{equation*}
  \hdim \mathfrak{B}(G,\mathfrak{p}) =\hdim \Lambda(G).
\end{equation*}
In addition, the exponent of convergence of $G$,
\begin{equation*}
  \delta(G):=  \inf\{s>0\colon \sum_{g\in 
    G}\lambda_g^{-s}<\infty\}=
  \hdim \Lambda(G)
\end{equation*}
\cite{BJ97,Patterson76b,Sullivan84}.  In view of
Lemma~\ref{lem:HCLemma}, this is perhaps not so surprising.  In a
striking parallel with continued fractions, badly approximable
points correspond to bounded orbits of flows on
manifolds~\cite{BJ97,Dani86a,DaniNTDS,FM95,Series85}.

\subsection{Jarník--Besicovitch theorem}
\label{sec:JBt}

When $v>1$, the set $\mathcal{D}_v$ is complementary to the set $R_v$
of points approximable to exponent $v$, which is related to the set
$W_v$ of $v$-approximable numbers (these statements also hold for the
higher dimensional analogues). The Hausdorff dimension of $W_v$ was
determined by Jarník in 1929~\cite{Ja29} and independently by
Besicovitch in 1934~\cite{Bes34}.
\begin{theorem}[Jarník--Besicovitch] 
  \label{thm:JBthm}
  When $v\geqslant 1$, 
  \begin{equation*}
    \hdim W_v = \dfrac{2}{v+1},
  \end{equation*}
  and when $v \leq 1$, $W_v = \mathbb{R}$.
\end{theorem}
Establishing the upper bound is not difficult, since $W_v$ is a limsup
set. There is no loss of generality in working with the more
convenient set $W_v\cap[0,1]$, as 
$W_v=\bigcup_{k\in\mathbb{Z}}(W_v\cap[0,1] +k)$.  Consider the limsup set
\begin{displaymath}
  W_v\cap[0,1]=\bigcap_{N=1}^\infty \bigcup_{q=N}^\infty 
  \bigcup_{p=0}^q
  B(p/q,q^{-v-1}),
\end{displaymath}
where $B(p/q,\varepsilon)=\{x\in[0,1]\colon
\abs{x-p/q}<\varepsilon\}$, so that
$\di(B(p/q,\varepsilon)) \leqslant 2\varepsilon$. Then 
\begin{displaymath}
  \sum_{q=1}^\infty \sum_{p=0}^q \di(B(p/q,q^{-v-1}))^s \leqslant
  2^s \sum_{q=1}^\infty  q^{1-s(v+1)}<\infty
\end{displaymath}
when $s>2/(v+1)$.  Hence by Lemma~\ref{lem:HCLemma} and the properties 
of Hausdorff dimension, $\hdim W_v\leqslant 2/(v+1)$ when $v>1$.  
When $v \leq 1$, the theorem follows from Dirichlet's theorem.

Establishing the correct lower bound is much harder. Jarník's lengthy
and complicated proof involved continued fractions and arithmetic
arguments.  Besicovitch's proof was simpler and more geometric and is
the basis of regular and ubiquitous systems which have turned out to
be very effective techniques in determining the Hausdorff dimension of
a variety of sets~\cite{BS,bbdbaker,DRV90a}. Indeed ubiquity can imply
Khintchine's theorem~\cite{BDV03}, so ideas developed for the study of
the the Hausdorff dimension of the null sets also contribute
substantially to our understanding of this theorem. The
Jarník--Besicovitch theorem has been extended considerably, to higher
dimensions, hyperbolic space and local fields.  Jarník himself proved
the Hausdorff measure analogue of Khintchine's theorem for
simultaneous Diophantine approximation~\cite{Ja31} and deduced that
the set of points in $\mathbb{R}^n$ satisfying
\begin{displaymath}
  \abs{\boldsymbol{\alpha} - \dfrac{\mathbf{p}}{q}} < q^{-v-1}
\end{displaymath}
for infinitely many $\mathbf{p} \in \mathbb{Z}^n$, $q \in \mathbb{N}$
has Hausdorff dimension $(n+1)/(v+1)$ when $v \geqslant 1/n$ and $n$
otherwise. He also showed that the Hausdorff $s$-measure of $W_v$ at
the critical dimension is infinite. In view of the
inclusions~\eqref{eq:UEWincl}, it follows that the same holds for
$E_v$, $R_v$.  
The points in $W_v$ form an uncountable totally disconnected subset of
the line and so are in Mandelbrot's picturesque language `fractal
dust', as is $\mathfrak{B}$. 

We have seen that the notion of a very well approximable point extends
naturally to systems of real and $p$-adic linear forms and to
hyperbolic space (where they can be interpreted in terms of geodesic
excursions~\cite{Dani85}).  The Jarník--Besicovitch theorem and the
Hausdorff measure analogue of Khintchine's theorem have been
established in the real case~\cite{HDSV97,mmd92}, the $p$-adic
case~\cite{AA95,ddj99}, the formal power series
case~\cite{kristensen02} and the hyperbolic case~\cite{HVJBGF,MP93}.
Other generalisations are to restricted sequences~\cite{Rynne92b},
inhomogeneous Diophantine approximation~\cite{JL98} and to small
linear forms~\cite{HD93}.  A further generalisation to `shrinking
targets' has revealed some unexpected connections with complex
dynamics and ergodic theory~\cite{HV95}.

We will be interested in the case of a single real linear form
in~\S\ref{sec:applns}.  The set of $v$-approximable points
$\boldsymbol{\alpha}\in\mathbb{R}^n$ (see \eqref{eq:valf}) will be
denoted $L_v$, \emph{i.e.},
\begin{equation}
  \label{eq:Lv}
  L_v=\{\boldsymbol{\alpha}\colon
  \abs{\mathbf{q}\cdot\boldsymbol{\alpha}-p}<\abs{\mathbf{q}}_\infty^{-v}
  \text{ for infinitely many } \ \mathbf{q}\in\mathbb{Z}^n,
  p\in\mathbb{Z}\,\}
\end{equation}
and by the dual or linear form version of the 
Jarník--Besicovitch theorem~\cite{BD86},
\begin{equation}
  \label{eq:dimLv}
  \hdim L_v=
  \begin{cases}
    n-1+\dfrac{n+1}{v+1} &\text{when } \ v>n, \\
    n                  &\text{when }\  v\leqslant n. \\
  \end{cases}
\end{equation}
The $n-1$ term arises from the dimension of the resonant hyperplanes
\begin{displaymath}
  \{\mathbf{x}\in\mathbb{R}^n\colon  \mathbf{q}\cdot \mathbf{x}=p\}  
\end{displaymath}
sets while the other represents `fractal dust' normal to the
hyperplanes.

\subsection{Approximation by ratios of Gaussian integers}
\label{sec:DAGI}
Recall that the Gaussian integers are defined as the set
$\mathbb{Z}[i] = \{p_1 + i p_2 \in \mathbb{C} : p_1, p_2 \in
\mathbb{Z}\}$. These form a ring (in fact a unique factorisation
domain). The Gaussian rationals are defined as the set $\mathbb{Q}(i)
= \{a/b + i c/d \colon a/b, c/d \in \mathbb{Q}\}$.
Approximation by Gaussian rationals decouples into the independent
approximation of the real and imaginary part respectively. Here we
study the more interesting problem of approximation by ratios of
Gaussian integers.

Approximation of complex numbers by ratios of Gaussian integers was
studied by Hermite and Hurwitz in the 19th
century~\cite[IV,~\S4]{Koksma} but, unlike in the real case, a
continued fraction approach did not give the best possible analogue of
Dirichlet's theorem.  This was obtained in 1925 by
Ford~\cite{Ford25}, who used additional geometrical ideas based
on the Picard group $\SL(2,\mathbb{Z}[i])$. The Gaussian analogues of
the theorems of Dirichlet, Khintchine, Jarn\'ik and Besicovitch
treated below should be compared with those corresponding to
simultaneous Diophantine approximation in the real plane
$\mathbb{R}^2$. 

In 1967 A.~Schmidt introduced a theory of regular and dually regular
chains for continued fractions, to treat approximation problems in
complex numbers~\cite{ASchmidt67,ASchmidt75a}.  Schmidt was concerned
with the study of complex quadratic irrationals and the complex
version of Pell's equation and with extensions to
groups~\cite{ASchmidt69,ASchmidt74}.  Our interest, however, is with
complex or Gaussian rational analogues of Dirichlet's theorem and with
the observation that the extended complex plane is the limit set of
the Picard group.  In this connection, Patterson established analogues
of Dirichlet's theorem for the less general Fuchsian
groups~\cite{Patterson76a,Velani93}.  Later Stratmann and Velani
obtained versions for Kleinian
groups~\cite{Stratmann94,SV95,VelaniKTHD} and so for the Picard group.
These results can be translated into complex versions of Dirichlet's
theorem with an undetermined constant.  Nevertheless, for
completeness, a short geometry of numbers proof of the complex version
of Dirichlet's theorem is given below.  Although the constant here is
not best possible, the result is all we need.  Proofs of the complex
analogues of Jarník's theorem on badly approximable numbers and the
Jarník--Besicovitch theorem that do not use the hyperbolic space
framework will be given in \S\ref{sec:cba} and \S\ref{sec:cJBt} below.

Complex Diophantine approximation has also been investigated from the
point of view of the distribution of the values of polynomials with
real integer coefficients but with complex variable $z$~\cite{bd99};
for another complex analogue see~\cite{Harman02}. \emph{For the rest
  of this section, $p=p_1+ip_2$, $q=q_1+iq_2$ will denote Gaussian
  integers with $q\ne 0$}.
\begin{theorem}
  Given any $z=x+iy\in \mathbb{C}$ and $N\in \mathbb{N}$, there exist
  Gaussian integers $p=p_1+ip_2$, $q=q_1+iq_2$ with $0<\vert q\vert
  \leqslant N$ such that
  \begin{equation}
    \label{eq:CDT2}
    \left\vert z-\dfrac{p}{q}\right\vert  < \dfrac{2}{\vert q\vert N}.
  \end{equation}
  Moreover for infinitely many $p, q\in \mathbb{Z}[i]$,
  \begin{equation}
    \label{eq:Dio}
    \left\vert z-\dfrac{p}{q}\right\vert  <  \dfrac{2}{\vert q\vert
      ^2}. 
  \end{equation}
\end{theorem}

\begin{proof} 
  The inequality~\eqref{eq:CDT2} holds if and only if the inequality
  \begin{equation}
    \label{eq:CDT}
    \left\vert x+iy-\dfrac{p_1+ip_2}{q_1+iq_2}\right\vert  
    < \dfrac{2}{\vert q_1+iq_2\vert N} 
  \end{equation}
  holds, \emph{i.e.}, if and only if
  \begin{equation*}
    \label{eq:CDT1}
    \left\vert (q_1x-q_2y-p_1)+i(q_2x+q_1y-p_2)\right\vert 
    <\dfrac{2}{N}
  \end{equation*}
  holds, which is the case if 
  \begin{equation}
    \label{eq:CDTht}
    \max\{\vert q_1x-q_2y-p_1\vert , \vert q_2x+q_1y-p_2\vert
    \}< \dfrac{\sqrt{2}}{N}. 
  \end{equation}
  By Minkowski's linear forms theorem, the system of inequalities
  \begin{eqnarray*}
    \vert q_1x-q_2y-p_1\vert  & < & 2^{1/2} \, N^{-1} \\
    \vert q_2x+q_1y-p_2\vert  & < &  2^{1/2} \, N^{-1} \\
    \vert q_1\vert  & \leqslant &  2^{-1/2} N \\
    \vert q_2\vert  & \leqslant &  2^{-1/2} N
  \end{eqnarray*}
  has a non-zero solution in integers $p_1,p_2,q_1,q_2$.
  Hence~\eqref{eq:CDTht} has a solution with $\vert q\vert =\vert
  q_1+iq_2\vert \leqslant N$, as claimed.
\end{proof}
Since the Gaussian rationals $p/q$ are not required to be on lowest
terms, \eqref{eq:Dio} holds infinitely often. For if $p/q$
fails to satisfy~\eqref{eq:Dio}, then \emph{a fortiori}, $(\kappa
p)/(\kappa q)$, where $\kappa$ is a non-zero Gaussian integer, will
also fail. And if $p/q$ satisfies~\eqref{eq:Dio}, only a finite number
of the fractions $(\kappa p)/(\kappa q)$ can also satisfy it.

This result should be compared to \eqref{eq:5}. As in the real case,
complex numbers for which Dirichlet's theorem cannot be significantly
improved are called badly approximable and numbers for which it can
are called very well approximable. More precisely, a complex number
$z$ is \emph{badly approximable} if there exists a constant $K=K(z)$
such that for all $p,q \in \mathbb{Z}[i]$, $q \neq 0$, 
\begin{equation*}
  \left\vert z-\dfrac{p}{ q}\right\vert \geqslant \dfrac{K}{\vert
    q\vert ^2}. 
\end{equation*}
Ford~\cite{Ford25} showed that the complex quadratic irrationals
$(1\pm i\sqrt3)/2$ are the worst approximable numbers with
$K=1/\sqrt3$ and thus correspond to the golden ratio $(\sqrt5+1)/2$ in
the real case.  Badly approximable numbers have been studied by
A.~Schmidt from the viewpoint of the Markoff
spectrum~\cite{ASchmidt84}.

As in the real case, a complex number $z$ is 
\emph{$v$-approximable} if
\begin{equation}
\label{eq:vapprox}
  \left\vert z-\dfrac{p}{q}\right\vert < \dfrac{1}{\vert q\vert ^{v+1}}
\end{equation}
for infinitely many $p, q \in \mathbb{Z}(i)$ and will be called
\emph{$\Psi$-approximable} if 
\begin{equation}
  \label{eq:ineq1}
  \left\vert z-\dfrac{p}{q}\right\vert 
  <\Psi(\vert q\vert),
\end{equation} 
where $\Psi: [1, \infty) \rightarrow \mathbb{R}^+$. The set of
$\Psi$-approximable $z$ will be written $W^*(\Psi)$. Thus $W^*(\Psi)$
is the set of points $z\in\mathbb{C}$ for which the
inequality~\eqref{eq:ineq1} holds for infinitely many Gaussian 
rationals $p/q$. The set of $v$-approximable complex numbers will be
written $W_v^*$.

\subsection{Khintchine's theorem for complex numbers}

The complex analogue of Khintchine's theorem was proved in 1952 by
LeVeque~\cite{leveque1952} who combined Khintchine's continued
fraction approach with ideas from hyperbolic geometry.  In 1976,
Patterson proved slightly less sharp versions of Khintchine's theorem
for Fuchsian groups acting on hyperbolic
space~\cite{Patterson76a,Velani93}.  A little later,
Sullivan~\cite{Sullivan82} used Bianchi groups and some powerful
hyperbolic geometry arguments to prove more general Khintchine
theorems for real and for complex numbers.  In the latter case, the
result includes approximation of complex numbers by ratios $a/b$ of
integers $a,b$ from the imaginary quadratic field
$\mathbb{R}(i\sqrt{d})$, where $d$ is a squarefree natural number. The
case $d=1$ corresponds to the Picard group and approximation by ratios
of Gaussian integers. Stratmann and Velani extended Patterson's
results with similar minor technical restrictions to Kleinian
groups~\cite{SV95,VelaniKTHD}. These include the Bianchi groups and
give a less precise and differently formulated version of Sullivan's
result. We now state LeVeque's result in our notation.
\begin{theorem}
  \label{thm:complexKT}
  Suppose $k^2\Psi^(k)$ is decreasing. Then the Lebesgue measure of
  $W^*(\Psi)$ is null or full according as the sum
  \begin{equation}
    \label{eq:Ksum}
    \sum_{k=1}^\infty k^3\Psi(k)^2
  \end{equation}
  converges or diverges.
\end{theorem}
As well as being more general, Sullivan's result is more precise as
the growth condition for the function $\Psi$ is weakened to a
`comparability condition'. Instead of the sum, Sullivan and LeVeque
use the equivalent integral $\int x^3\Psi(x)^2\,dx$ (their integrands
are different owing to different forms of~\eqref{eq:ineq1}).  It is
readily verified that $W^*(\Psi)$ is invariant under translations by
Gaussian integers $p=p_1+ip_2$, so that
\begin{equation}
  \label{eq:2}
  W^*(\Psi)= \bigcup_{p\in\mathbb{Z}[i]} V^*(\Psi)+p = 
  \bigcup_{p_1,p_2\in\mathbb{Z}}
  V^*(\Psi)+p_1+ip_2,
\end{equation}
where $V^*(\Psi)=W^*(\Psi)\cap I^2$ and $I^2=[0,1)^2=\{x+iy\colon
0\leqslant x,y < 1 \}$.  We will work in the more convenient unit
square $I^2$ and consider the set $V^*(\Psi)$.

Before continuing, we need some more definitions and notation.

\subsection{Resonant sets  and balls in $I^2$}

Let $q=q_1+iq_2\in\mathbb{Z}[i]\setminus\{0\}$. The set
$R(q)=R(q_1,q_2)\subset I^2$ where
\begin{align*}
  \label{eq:resetal}
  R(q)&=\left\{\dfrac{p}{q}\colon p\in\mathbb{Z}[i]\right\}\cap I^2
  =
  R(q_1,q_2)
  \\&=
  \left\{\left(\dfrac{p_1q_1+p_2q_2}{q_1^2+q_2^2}, 
      \dfrac{p_2q_1-p_1q_2}{q_1^2+q_2^2}\right)
    \colon p_1,p_2\in\mathbb{Z}\right\}\cap I^2
\end{align*}
is called a \emph{resonant set}. This set is the analogue in the
complex plane of the set $\{p_1/q_1\colon 0\leqslant p_1<q_1\}$ in the
real line.  The points in the resonant set form a lattice inclined at
an angle $\tan^{-1} (q_1/q_2)$ to the real axis and in which the side
length of the fundamental region is $\vert q\vert
=(q_1^2+q_2^2)^{-1/2}$.  Area and congruence considerations give that
the number of points of $ R(q_1,q_2)$ in $I^2$ is
\begin{equation}
  \label{eq:pts}
  \# R(q)= \# R(q_1,q_2)=\vert q\vert ^2 
  = q_1^2+q_2^2.  
\end{equation}

\begin{figure}[htbp]
  \centering
  \includegraphics[width=6cm,keepaspectratio]{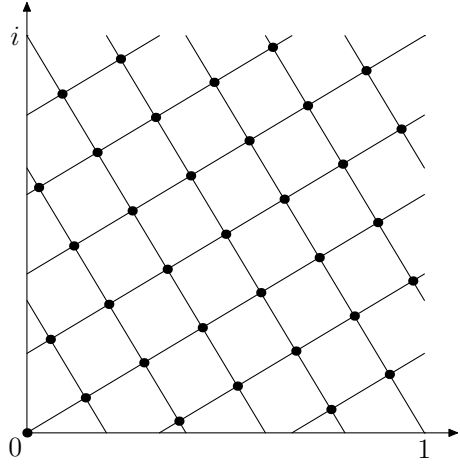}
  \caption{The lattice in $I^2$ corresponding to $q_1=5$, $q_2=3$.} 
  \label{fig:latt}
\end{figure}

The disc 
\begin{equation}
  \label{eq:disc}
  D(p/q;\varepsilon)
  =\left\{z\in \mathbb{C}\colon
    \left\vert z-p/q\right\vert < \varepsilon
  \right\},
\end{equation}
with radius $\varepsilon$ and centred at $p/q$ where
\begin{equation*}
  \dfrac{p}{q}=\dfrac{p_1+ip_2}{q_1+iq_2}=
  \dfrac{(p_1+ip_2)(q_1-iq_2)}{q_1^2+q_2^2}=
  \dfrac{p_1q_1+p_2q_2}{q_1^2+q_2^2} +
  i\dfrac{p_2q_1-p_1q_2}{q_1^2+q_2^2} 
\end{equation*}
has area $\pi \varepsilon^2$. The set
\begin{equation}
  \label{eq:Bcd}
  B(q,\varepsilon)=\left\{z\in I^2 \colon
    \left\vert z-\dfrac{p}{q}\right\vert <\varepsilon \text{ for some } p\in
    \mathbb{Z}[i]\right\}=\bigcup_{p\in\mathbb{Z}[i]}
  D(p/q,\varepsilon)
\end{equation}
can be regarded as a neighbourhood of a resonant set and its measure
\begin{equation}
  \label{eq:|Bcd|}
  \vert  B(q,\varepsilon)\vert = \pi \varepsilon^2 \left(\vert q\vert ^2+O(\vert q\vert )\right). 
\end{equation}
Now $V^*(\Psi)$ can be expressed in the form
\begin{equation}
  \label{eq:limsupW}
  V^*(\Psi)= \bigcap_{N=1}^\infty\bigcup_{k=N}^\infty 
  \bigcup_{\vert q\vert =k}
  B(q,\Psi(\vert q\vert)) =\limsup_{\vert q\vert \to\infty}
  B(q,\Psi(\vert q\vert)).
\end{equation}
Thus $V^*(\Psi)$ is a limsup set and for each $N=1,2,\dots$  has a 
natural cover
\begin{equation}
  \label{eq:cover}
  \mathcal{C}_N(V^*(\Psi))=  
  \{ B(q,\Psi(\vert q\vert)\colon \vert q\vert \geqslant N\}.
\end{equation}
Hence for each $N=1,2,\dots$, the measure of $V^*(\Psi)$ satisfies
\begin{displaymath}
  \vert V^*(\Psi) \vert \leqslant \sum_{k=N}^\infty \sum_{k \leqslant
    \abs{q} < k+1} \vert B(q,\Psi(\vert q\vert)\vert \ll
  \sum_{k=N}^\infty k^2 \Psi(k)^2 \sum_{k \leqslant \abs{q} < k+1} 1. 
\end{displaymath}
Now, $\sum_{1 \leqslant \abs{q} \leqslant k} 1$ is the number of
lattice points in the closed disc $D(0,k)$ of radius $k$. By
\cite{MR27:4799}, given any $\varepsilon > 0$, 
\begin{displaymath}
  \sum_{1 \leqslant \abs{q} \leqslant k} 1 = \pi k^2 + O\left(k^{12/37
      + \varepsilon}\right).
\end{displaymath}
Hence 
\begin{equation}
  \label{eq:7}
  \sum_{k \leqslant \abs{q} < k+1} 1 = \pi (k+1)^2 +
  O\left(k^{1/3}\right) - \pi k^2 + O\left(k^{1/3}\right) = 2 \pi k
  + O\left(k^{1/3}\right).
\end{equation}
Therefore
\begin{displaymath}
  \abs{V^*(\Psi)} \ll \sum_{k=N}^\infty k^3\Psi(k)^2.
\end{displaymath}
Since $N$ is arbitrary, the  convergence of the series $\sum_k
k^3\Psi(k)^2$ implies that
\begin{displaymath}
  \abs{V^*(\Psi)}=\abs{W^*(\Psi)}=0,
\end{displaymath}
which is the convergence part of the complex analogue of Khintchine's
theorem.  The much more difficult case of divergence requires deeper
arguments and the reader is referred to Sullivan's bold and highly
geometrical paper~\cite{Sullivan82}.

\section{Badly approximable complex numbers and Jarník's theorem}
\label{sec:cba}
The results of~\cite{BJ97} and of~\cite{FM95} on badly approximable
points arising with the action of Kleinian groups and on bounded
geodesics on Riemann surfaces respectively, specialised to the case of
the Picard group, could be translated to give Jarník's theorem for
complex numbers badly approximable by ratios of Gaussian integers.
However, we will give a self-contained proof using an extension to the
complex numbers of the $(\alpha,\beta)$-game introduced by
W.~M.~Schmidt~\cite{Schmidt66}. In the setting of the complex plane,
$S \subset \mathbb{C}$, $\alpha, \beta \in (0,1)$ and discs replace
intervals. Thus, the game begins with player B choosing a disc $B_1 =
\{ z \in \mathbb{C} : \abs{z-b_1} \leqslant \rho_1\}$. Next A chooses
a disc $A_1 \subset B_1$ with radius $\alpha \rho_1$. Then B chooses a
disc $B_2 \subset A_2$ with radius $\beta\alpha\rho_1$ and so on ad
infinitum, such that $B_{n+1}$ has radius $(\alpha \beta)^n \rho_1$
for any $n \geq 0$.

The discs $B_1, A_1, B_2, A_2, \dots$ form a nested decreasing
sequence of closed sets, so there is a unique intersection point,
$\omega$ say. Player A wins if this point is an element of the set
$S$, \emph{i.e.}, if $\bigcap_j A_j = \{\omega\} \subset S$.
Otherwise, B wins. A set $S$ is said to be $(\alpha, \beta)$-winning
if A can win the game for the parameters $\alpha$ and $\beta$ no
matter how well B plays. If for some $\alpha \in (0,1)$, A can win the
game for any $\beta \in (0,1)$, the set $S$ is said to be
$\alpha$-winning.

Player A benefits from $\alpha$ being small. Indeed, when $\alpha$
gets smaller, A can limit the amount of choice B has in the next
move. As in the real case, it may be shown that if $\alpha' < \alpha$
and $S$ is $\alpha$-winning, then $S$ is $\alpha'$-winning.  Hence,
given $S$, there is a largest value of $\alpha$ for which $S$ is
$\alpha$-winning so we may define
\begin{displaymath}
  \alpha^*(S) := \sup \left\{\alpha \in (0,1) : S \text{ is
      $\alpha$-winning} \right\}.
\end{displaymath}

In fact, for any $S \neq \mathbb{C}$, we easily see that $\alpha^*(S)
\leqslant \tfrac{1}{2}$. To see this, note that for the parameters
$\alpha > \tfrac{1}{2}$ and $\beta \in (0, 2\alpha -1)$, B may ensure
that the centres of all the $B_i$ are the same. Hence, by choosing the
first disc $B_1$ with centre $b_1 \notin S$, we have the result.

When $S = \mathfrak{B}$, A wins if she can force $\omega$ to be badly
approximable, \emph{i.e.}, if $\omega$ is in the set
\begin{displaymath}
  \mathfrak{B} = \left\{ z \in \mathbb{C} : \exists K > 0 \: \forall 
    p,q \in
    \mathbb{Z}[i]: \abs{z - \dfrac{p}{q}} > 
    \dfrac{K}{\abs{q}^2}\right\}.
\end{displaymath}
Since $\mathfrak{B}$ is null, this seems unlikely to be the case, but in
fact A can win the game whenever $2\alpha < 1 + \alpha\beta$.  It
immediately follows that
\begin{equation}
  \label{seq:5}
  \alpha^*(\mathfrak{B}) = 1/2,
\end{equation}
so in fact A may win this game almost as easily as she could win the
game where $S$ is the entire complex plane with one point removed. We
take some time to prove \eqref{seq:5}.

We have already seen that $\alpha^*(\mathfrak{B}) \leqslant 1/2$, so
we will only consider $\alpha \leqslant 1/2$, as this will simplify
the proof. Let $\alpha \in (0,1/2]$ and $\beta \in (0,1)$ be
fixed. Note that $\gamma := 1 + \alpha \beta - 2 \alpha > 0$. We may
assume without loss of generality that the radius of the initial ball
$B_1$ satisfies $\rho_1 \leqslant
\alpha\beta\tfrac{\gamma}{8}$. Indeed, otherwise we would let the game
continue in an arbitrary fashion until reaching a $B_j$ for which
$\rho_j \leqslant \alpha\beta\tfrac{\gamma}{8}$ and then take this to
be our starting point. The constant $K$ in the definition of
$\mathfrak{B}$ will be $\delta = \tfrac{\gamma}{4} \min(\rho, \alpha^2
\beta^2 \tfrac{\gamma}{8})$, where $\rho = \rho_1$.

Let $t \in \mathbb{N}$ be such that $\alpha \beta \gamma \leqslant 2
(\alpha\beta)^t < \gamma $. Let $R = (\alpha \beta)^{-t/2}$. Clearly,
it suffices to prove for any $n \in \mathbb{N}$ that if
\begin{subequations}
  \begin{eqnarray}
    \label{seq:3}
    \gcd(p,q) = 1, \\ z \in B_{nt+1}, \\
    \label{seq:4}
    0 < \abs{q} < R^n,
  \end{eqnarray}
\end{subequations}
then $\abs{z - \tfrac{p}{q}} > \tfrac{\delta}{\abs{q}^2}$.  This may
be done using induction.

For $n = 0$, there is nothing to prove, as \eqref{seq:4} leaves us no
$q$ to consider. Hence, we may assume that we have $B_1, \dots,
B_{(k-1)t+1}$ such that the above holds for $0 \leqslant n \leqslant
k-1$. In subsequent play, A thus only needs to worry about
$\tfrac{p}{q}$ with $R^{k-1} \leqslant \abs{q} < R^k$, as the
remaining problematic fractions have been sorted out in the preceding
steps of the game.

In fact, there can be at most one such $\tfrac{p}{q}$. Indeed, suppose
that there exists $p,p',q,q' \in \mathbb{Z}[i]$ and $z, z' \in B_{(k-1)t
  + 1}$ with
\begin{displaymath}
  \abs{z - \dfrac{p}{q}} \leqslant \dfrac{\delta}{\abs{q}^2} \quad
  \text{and} \quad \abs{z' - \dfrac{p'}{q'}} \leqslant 
  \dfrac{\delta}{\abs{q'}^2}
\end{displaymath}
with $\abs{q}, \abs{q'} \in [R^{k-1}, R^k)$. Then
\begin{multline*}
  \abs{\dfrac{p}{q} - \dfrac{p'}{q'}} \leqslant \abs{z- 
    \dfrac{p}{q}} +
  \abs{z' - \dfrac{p'}{q'}} + \abs{z-z'} \\ \leqslant
  \dfrac{\delta}{\abs{q}^2} + \dfrac{\delta}{\abs{q'}^2} +
  2\rho(\alpha\beta)^{(k-1)t} \leqslant 2 \delta R^{2-2k} + 2 \rho 
  R^{2-2k}
  \\ \leqslant 4 \rho (\alpha\beta)^{-t} R^{-2k} \leqslant 4 \tfrac{1}{8} 
  \alpha
  \beta \gamma \dfrac{2}{\alpha \beta \gamma} R^{-2k} = 
  R^{-2k}.
\end{multline*}
On the other hand, since $\gcd(p,q) = \gcd(p',q') = 1$ and since
$\mathbb{Z}[i]$ is a unique factorisation domain,
\begin{displaymath}
  \abs{\dfrac{p}{q} - \dfrac{p'}{q'}} = \abs{\dfrac{pq' - 
      p'q}{q q'}}
  \geqslant \dfrac{1}{\abs{q}\abs{q'}} > R^{-2k}
\end{displaymath}
whenever $p\neq p'$ and $q \neq q'$. Hence, there can be at most one
problematic point.

Note that we have used the property of the ring of Gaussian integers
being a unique factorisation domain. For other rings, this may not be
the case and stronger tools are needed for proving the analogous
result. However, for clarity of this exposition, we will take the
simple route and use the unique factorisation property. Note also
that for the corresponding result in simultaneous Diophantine
approximation, the module $\mathbb{Z}^2$ takes the place of the
Gaussian integers. Even though the underlying sets are the same, the
difference in algebraic structure prevents the method of this proof  
from working in the case of simultaneous Diophantine approximation.

As there can be at most one point in $B_{(k-1)t+1}$ suitably close to
a Gaussian rational $p/q$, we may devise strategies that avoids a disc
around this point $p/q$ of a suitable radius. It is clear that we need
to avoid $C=B(\tfrac{p}{q}, \tfrac{\delta}{\abs{q}^2})$. We examine
two possibilities in turn.

\begin{figure}[htbp]
  \centering 
  \includegraphics[width=6cm,keepaspectratio]{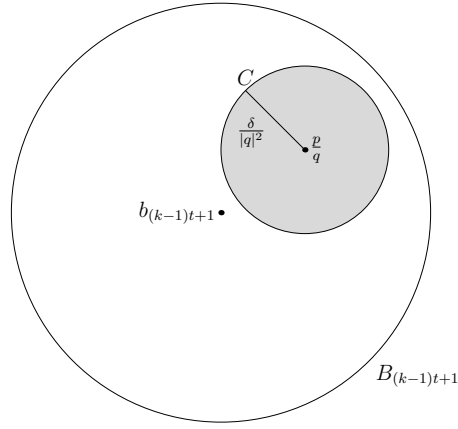}
  \caption{A avoids $C$ in one move}
  \label{fig:easywin}
\end{figure}
First, consider the case when $\abs{\tfrac{p}{q} - b_{(k-1)t+1}} >
\delta R^{2-2k}$. In this case,
\begin{displaymath}
  \abs{\dfrac{p}{q} - b_{(k-1)t+1}} > \delta R^{2-2k} \geqslant
  \dfrac{\delta}{\abs{q}^2}. 
\end{displaymath}
Hence, $b_{(k-1)t+1} \notin C$, so we are in the situation of figure
\ref{fig:easywin}. As $\alpha \leqslant 1/2$, we see that A can choose her
next disc in such a way that it does not intersect with $C$.

\begin{figure}[htbp]
  \centering
  \includegraphics[width=6cm,keepaspectratio]{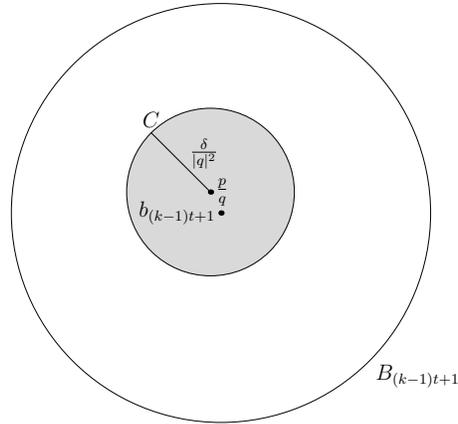}
  \caption{A avoids $C$ in $t$ steps}
  \label{fig:difficult}
\end{figure}

This leaves the final case when
\begin{displaymath}
  \abs{\dfrac{p}{q} - b_{(k-1)t+1}}\leqslant \delta R^{2-2k}.  
\end{displaymath}
This corresponds to the case when $b_{(k-1)t+1} \in C$ (figure
\ref{fig:difficult}).
In this case, we clearly need to work harder to get an answer, as the
right strategy of A is not immediately obvious. However, it turns out
that picking a fixed direction and moving as far as possible in this
direction at each turn will cause the disc chosen by B after $t$ steps
to have empty intersection with $C$. The following lemma formalises
this.

\begin{lemma}
  \label{lem:fundamental}
  Let $\alpha, \beta \in (0,1)$ with $\gamma = 1 + \alpha\beta -
  2\alpha >0$. Let $t \in \mathbb{N}$ be such that $(\alpha\beta)^t <
  \gamma/2$. Suppose that the disc $B_k = (b_k, \rho_k)$ occurs at
  some stage in an $(\alpha, \beta)$-game. Then A can play in such a
  way that
  \begin{displaymath}
    B_{k+t} \subseteq \left\{z \in \mathbb{C} : \abs{z-b_k} > 
      \rho_k \tfrac{\gamma}{2} \right\}.
  \end{displaymath}
\end{lemma}

\begin{proof}
  We define a strategy for A in the following way: Suppose that the
  last disc chosen by B was $B(c,\rho)$ for some $\rho > 0$. Choose
  some $\hat{c} \in \mathbb{C}$ with  $\abs{\hat{c}} =
  \rho(1-\alpha)$. We define a legal move for A by the map
  \begin{displaymath}
    B(c,\rho) \mapsto B(c+\hat{c}, \alpha \rho).
  \end{displaymath}
  We need to convince ourselves that this move is in fact a legal
  move. First, we note that the radius is the right one. Hence, we
  need only prove the inclusion $B(c+\hat{c}, \alpha \rho) \subseteq
  B(c,\rho)$. But this follows since for any $z \in B(c+\hat{c},
  \alpha \rho)$,
  \begin{displaymath}
    \abs{z - c} = \abs{z - \hat{c} + \hat{c} - c} \leqslant \abs{z -
      (c + \hat{c})} + \abs{\hat{c}} \leqslant \alpha \rho + \rho
    (1-\alpha) = \rho.
  \end{displaymath}
  The strategy of player A will be to use the above move, no matter
  how B plays the game. 
  
  We denote the discs chosen by B by $B_k = (b_k, \rho_k)$, and the
  discs chosen by A by $A_k = (a_k, \alpha \rho_k)$. Note that
  \begin{equation}
    \label{seq:1}
    \abs{a_k - b_k} = \abs{b_k + \hat{c} - b_k} = 
    \rho_k(1-\alpha).
  \end{equation}
  \begin{figure}[htbp]
    \centering 
    \includegraphics[width=6cm,keepaspectratio]{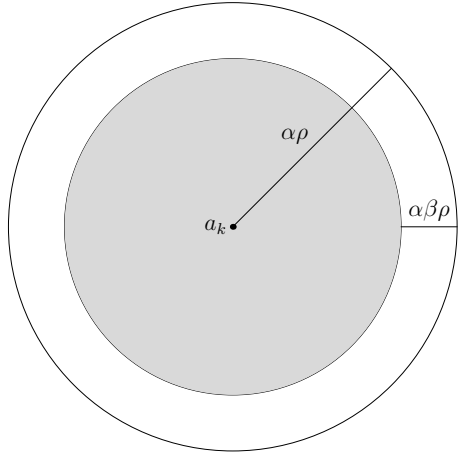}
    \caption{The point $b_{k+1}$ must be chosen in the shaded 
      area.}
    \label{fig:Blacks_choice}
  \end{figure}  
  Also, from figure \ref{fig:Blacks_choice}, we see that
  \begin{equation}
    \label{seq:2}
    \abs{b_{k+1} - a_k} \leqslant \alpha \rho_k - \alpha \beta \rho_k.
  \end{equation}
  Hence, by \eqref{seq:1} and \eqref{seq:2},
  \begin{multline*}
    \abs{b_k - b_{k+1}} = \abs{(b_k - a_k) - (b_{k+1} - a_k)}
    \geqslant \big\vert\abs{b_k - a_k} - \abs{b_{k+1} - a_k}
    \big\vert\\ 
    = \big\vert \rho_k (1-\alpha) - \abs{b_{k+1} - a_k}\big\vert >
    \abs{\rho_k (1-\alpha) - \alpha \rho_k + \alpha \beta \rho_k} =
    \rho_k \gamma. 
  \end{multline*}
  Continuing as above, at each step choosing $\hat{c}$ to point in the
  direction of the first one chosen, we obtain,
  \begin{displaymath}
    \abs{b_{k+t} - b_k} > \rho_k \gamma.
  \end{displaymath}
  But since $\rho_{k+t} = (\alpha \beta)^t \rho_k < \rho_k
  \tfrac{\gamma}{2}$, we have for any $z \in B_{k+t}$,
  \begin{displaymath}
    \abs{z - b_k} \geqslant \big\vert\abs{z - b_{k+t}} - \abs{b_{k+t}
      -  b_k}\big\vert > \rho_k \tfrac{\gamma}{2}.
  \end{displaymath}
  This completes the proof.
\end{proof}

With the induction step done, we have shown that A can play in such a
way that \eqref{seq:5} holds for the final intersection point.  In
fact, we have found an explicit lower bound on the constant $K$ for
which \eqref{seq:5} holds.

From \eqref{seq:5}, we can get the Hausdorff dimension of the set
$\mathfrak{B}$. By considering the amount of choice B has in the game
and using this to find a lower bound on the $s$-length of appropriate
covers, we obtain a lower bound on the dimension of any $(\alpha,
\beta)$-winning set. This construction was carried out in considerable
generality by W. M. Schmidt \cite{Schmidt66}, who obtained a lower
bound for the Hausdorff dimension of $(\alpha,\beta)$-winning sets in
real Hilbert spaces. We sketch an approach to obtaining the dimension
from the number $\alpha^*(\mathfrak{B})$. Because of the geometric
nature of the problem, $\mathfrak{B}$ will be regarded as a subset in
$\mathbb{R}^2$. We will compute the usual planar Hausdorff dimension
here and subsequently.

To obtain the Hausdorff dimension of $\mathfrak{B}$, we consider the
game from B's point of view. Assume first without loss of generality
that $B_0$, the first disc chosen, has radius $1$. At any point of the
game, B can choose to direct the game into a number of disjoint
discs. While there may be a variety of ways in which these discs can
be chosen, the maximum number $N(\beta)$ is roughly equal to
$1/\beta^2$, \emph{i.e.}, $N(\beta) \asymp 1/\beta^2$.

We limit B's choice to these $N(\beta)$ possible moves and assume that
A is playing to win the game. This gives us a parametrisation of the
sequence of discs chosen by B, so that $B_k=B_k(j_1, \dots, j_k)$ with
the $j_i \in \{0, \dots, N(\beta)-1\}$ for $i = 1, \dots, k$. For
later use, note that the radius of each these discs is $\rho_k =
(\alpha\beta)^k$. By simultaneously considering all the different
ways, B may play the game, we obtain a function
\begin{displaymath}
  f: \left\{0, \dots, N(\beta)-1\right\}^{\mathbb{N}} \rightarrow
  \mathfrak{B}, \quad (\lambda_k)_{k \in \mathbb{N}} \mapsto
  \bigcap_{k \in \mathbb{N}} B_k(\lambda_1, \dots, \lambda_k) =
  \{x(\lambda)\}.
\end{displaymath}
We define the set $\mathfrak{B}^* \subseteq \mathfrak{B}$ to be the
range of $f$. As every number in the interval $[0,1]$ has at least one
expansion in base $N(\beta)$, we may map this set onto the unit
interval by the map
\begin{displaymath}
  g: \mathfrak{B}^* \rightarrow [0,1], \quad x(\lambda) \mapsto
  0.\lambda_1 \lambda_2 \dots.
\end{displaymath}
Note that these functions could well be multivalued, but this is of no
concern to us. All we need is a cover of $[0,1]$. We extend this
function to subsets of the complex plane by defining $g(Z) = g(Z \cap
\mathfrak{B}^*)$ for any $Z \subseteq \mathbb{C}$, where by convention
$B(\emptyset) = 0$.

Now, take some cover $\mathcal{C} = \{C_l\}_{l \in \mathbb{N}}$ of
$\mathfrak{B}$ with discs of radius $\rho(C_l) = \rho_l$. We wish to
find a lower bound on the $s$-length of this cover for appropriate $s$
as $\rho_l$ becomes smaller. This is where the function $g$ comes in
handy. As $\mathcal{C}$ covers $\mathfrak{B}^*$, we see that
$g(\mathcal{C})$ covers $[0,1]$. We let $\overline{\mu}$ denote the
outer Lebesgue measure. By sub-additivity,
\begin{equation}
  \label{seq:10}
  \sum_{l=1}^\infty \overline{\mu} (g(C_l)) \geqslant
  \overline{\mu}\left(\bigcup_{l=1}^\infty g(C_l)\right) \geqslant 1.
\end{equation}

Now let $\omega>0$ be sufficiently small so that any disc of radius
$\omega(\alpha\beta)^k$ intersects at most two of the discs $B_k(j_1,
\dots, j_k)$. By \cite[Lemma 20]{Schmidt66}, $\omega = 2/\sqrt{3} -1 <
1$ has this property. We define integers
\begin{displaymath}
  k_l = \left[\dfrac{\log \left(2 \omega^{-1}
        \rho_l\right)}{\log(\alpha\beta)}\right].
\end{displaymath}

For $\rho_l$ sufficiently small, we see that $k_l > 0$ and $\rho_l <
\omega (\alpha\beta)^{k_l}$. Hence, the disc $C_{l}$ intersects at
most two of the discs $B_{k_l}(j_1, \dots, j_{k_l})$. As
$g(B_{k_l}(j_1, \dots, j_{k_l}))$ is clearly an interval of length
$N(\beta)^{-k_l}$, we have $\overline{\mu} (g(C_{l})) \leqslant 2
N(\beta)^{k_l}$, so by \eqref{seq:10},
\begin{displaymath}
  1 \leqslant \sum_{l=1}^\infty \overline{\mu} (g(C_l)) \leqslant
  \sum_{l=1}^\infty 2 N(\beta)^{k_l} \leqslant
  2(2\omega^{-1})^{\tfrac{\log(N(\beta))}{\abs{\log(\alpha\beta)}}}
  \sum_{l=1}^\infty
  \rho_l^{\tfrac{\log(N(\beta))}{\abs{\log(\alpha\beta)}}}.
\end{displaymath}
Thus, for $s = \log(N(\beta))/\abs{\log(\alpha\beta)}$, the $s$-length
of the cover $\mathcal{C}$ is strictly positive, so the Hausdorff
dimension of $\mathfrak{B}^*$ must be greater than this number.

Now, we fix $\alpha \in (0,1/2)$ and apply the above,
\begin{displaymath}
  \hdim(\mathfrak{B}) \geqslant \dfrac{\log \beta^{-2}}{\abs{\log
      \alpha \beta}} = \dfrac{2\abs{\log \beta}}{\abs{\log
      \alpha} + \abs{\log \beta}} \rightarrow 
      2 
\end{displaymath}
as $\beta \to 0$. We have thus proved the analogue of Jarník's theorem
for the complex numbers:
\begin{theorem}
  \label{thm:jarnik}
  The set $\mathfrak{B}$ is thick, \emph{i.e.}, for any disc $B \subseteq
  \mathbb{C}$,
  \begin{displaymath}
    \hdim(B \cap \mathfrak{B}) = 2.
  \end{displaymath}
\end{theorem}

\section{The complex Jarník--Besicovitch theorem}
\label{sec:cJBt} 
Let $W^*_v$ be the set of $v$-approximable complex numbers which
satisfy the inequality~\eqref{eq:vapprox} for infinitely many Gaussian
rationals $p/q$ (recall that $p,q \in \mathbb{Z}[i]$).  The general
Jarn\'ik--Besicovitch theorems in~\cite{HVJBGF,MP93} can be
specialised to  the Picard group to yield the Hausdorff dimension of
$W^*_v$ (see~\cite[Corollary~2]{MP93}) but we give a more
self-contained and direct proof. For convenience we consider
\begin{equation*}
  \label{eq:Vv}
  V^*_v:=W^*_v\cap I^2=\left\{z\in I^2\colon \left|z-\frac{p}{q}\right|
    <\frac1{|q|^{v+1}} \ {\text{ for
        infinitely many }} \frac{p}{q} \right\}.
\end{equation*}
\begin{theorem}
  \label{thm:complexJBt}
  \begin{equation*}
    \hdim V^*_v=\hdim W^*_v =
    \begin{cases}
      \dfrac{4}{v+1} & \text{ when } v\geqslant 1, \\
      2  & \text{ when } v\leqslant 1.
    \end{cases}
  \end{equation*}
\end{theorem}
When $v\leqslant 1$, $V^*_v$ is full by the complex form of
Khintchine's theorem so the second equality holds by (iii) in
\S~\ref{sec:HD}. As usual, the proof for $v > 1$ falls into two
parts. First, to obtain the upper bound for $\hdim V^*_v$, consider
the cover ${\scr{C}}(V^*_v)$ given by~\eqref{eq:cover} with $N=1$ and
$\Psi(k)=k^{-v-1}$. By \eqref{eq:7}, this has $s$-length
\begin{eqnarray*}
  \ell^s(\scr{C}(V^*_v)) & \ll & \sum_{k=1}^\infty \sum_{k \leqslant 
    \abs{q} < k+1} \abs{q}^2 \left(\abs{q}^{-v-1}\right)^s \ll
  \sum_{k=1}^\infty k^{2-(v+1)s} \sum_{k \leq \abs{q} < k+1} 1
  \\
  &\ll& \sum_{k=1}^\infty k^{3-s(v+1)} <\infty
\end{eqnarray*}
for $s>4/(v+1)$. It follows that when $v\geqslant 1$, 
\begin{equation}
  \label{eq:dimle}
  \hdim V^*_v\leqslant \dfrac{4}{v+1}.  
\end{equation}

To obtain the lower bound in the Jarník--Besicovitch theorem, we use
ubiquity \cite{MDAM, DRV90a}. Let $S\subseteq I^2$ and let $\rho\colon
\mathbb{N}\to (0,\infty)$ be a function.  Put
\begin{equation*}
  B(S;q,\varepsilon):=  
  \left\{z\in S \colon \left\vert z-\dfrac{p}{q}\right\vert <\varepsilon
    \ \text{ for some } \ p \in \mathbb{Z}[i]\right\}.
\end{equation*}
The set 
\begin{equation*}
  \mathcal{R}=\bigcup_{q} R(q)\subset \mathbb{C} ,  
\end{equation*}
where the union is over non-zero $q\in \mathbb{Z}[i]$, consists of
discrete points and so has dimension 0.  Let $S$ be any open square in
$I^2$.  By definition, $\mathcal{R}$ is \emph{ubiquitous} in
$S\subseteq I^2 =[0,1)^2$ with respect to $\rho$ if
\begin{equation*}
  \left\vert \bigcup_{q} B(S;q,\rho(N))\right\vert  \to \vert S\vert 
\end{equation*}
as $N\to \infty$. Now by the complex analogue of Dirichlet's theorem,
for each $N\in\mathbb{N}$,
\begin{eqnarray*}
  \left\{z\in S \colon \left\vert z-\dfrac{p}{q}\right\vert
    <\dfrac{2}{\vert q\vert N}\ \text{ for some }
    p,q\in\mathbb{Z}[i],1\leqslant\vert q\vert \leqslant N
  \right\} \\
 \qquad \qquad  =\bigcup_{q} B\left(q,\dfrac{2}{\vert q\vert N
    }\right) = S,\qquad\qquad\qquad\qquad\qquad &
\end{eqnarray*}
so that 
\begin{equation*}
  \left\vert \bigcup_{q} B\left(q,\dfrac{2}{\vert q\vert N
      }\right)\right\vert =\vert S\vert.
\end{equation*}
Consider the set $S(N)$ of $z\in S$ with `small denominators':
\begin{alignat*}{2}
  S(N) &= \left\{z\in S\colon \text{ there exist } p, q \text{ such 
      that } \left\vert z-\dfrac{p}{q}\right\vert <\dfrac{2}{\vert
      q\vert N}, 1\leqslant \vert q\vert <\dfrac{N}{\log N}\right\} \\ 
  &= \bigcup_{1\leqslant \vert q\vert < N/\log 
    N} B\left(q,\dfrac{2}{\vert q\vert N}\right).
\end{alignat*}
The measure $\vert S(N)\vert $ of $S(N)$ satisfies
\begin{eqnarray*}
  \abs{S(N)} &=&\sum_{1\leqslant \vert q\vert < N/\log N} \vert
  B(q,2/\vert q\vert N)\vert  \ll \sum_{1\leqslant \vert q\vert <
    N/\log N} (\vert q\vert N)^{-2}\vert q\vert ^2\\ 
  &\ll & N^{-2}\sum_{1\leqslant \vert q\vert < N/\log N} 1\ll N^{-2} 
  \left(N/\log N\right)^2 \ll (\log N)^{-2}\to 0
\end{eqnarray*}
as $N\to \infty$.  Choose $\rho(N)=2N^{-2}\log N$.  Then since
$\vert q\vert >N/\log N$ implies $\rho(N)>2/\vert q\vert N$,
\begin{equation*}
  \bigcup_{N/\log N < \vert q\vert \leqslant N} B(q,\rho(N))\supset 
  \bigcup_{N/\log N < \vert q\vert \leqslant N} B\left(q,\dfrac{2}{\vert
    q\vert N}\right) \to S 
\end{equation*}
in measure as $N\to\infty$ and so $\mathcal{R}$ is ubiquitous with
respect to $\rho(N)= 2\log N /N^2$ for any $S\subseteq I^2$ and so for
$S=I^2$. But since $\Psi$ is decreasing, by~\cite{DRV90a},
\begin{equation}
  \label{eq:dimge}
  \hdim V^*_v \geqslant \dim \mathcal{R} + \codim \mathcal{R}
  \limsup_N \dfrac{\log \rho(N)}{\log \Psi (N)} = \dfrac{4}{v+1}, 
\end{equation}
where $\dim$ is the topological dimension, so that $\dim \mathcal{R}
=0$, the codimension $\codim \mathcal{R}$ in $\mathbb{C}$ (regarded as
$\mathbb{R}^2$) of $\mathcal{R}$ is $2$ and $\Psi(N)=N^{-v-1}$.
                   
The required result follows on combining the two complementary
inequalities~\eqref{eq:dimle} and~\eqref{eq:dimge}.

In fact, since $S$ was an arbitrary open square, $\mathcal{R}$ is
\emph{locally} ubiquitous and it has been shown that the local
ubiquity of $\mathcal{R}$ also implies the divergence case of
Khintchine's theorem~\cite{BDV03}.

Note that even though one might expect Diophantine approximation in
$\mathbb{C}$ with respect to Gaussian integers to be similar to
simultaneous Diophantine approximation in the real plane
$\mathbb{R}^2$, this is not at all the case. The analogues of
Dirichlet's theorem, Khintchine's theorem and the
Jarn\'ik--Besicovitch theorem are quite different in the two cases.
Indeed the complex Dirichlet's Theorem and the Jarn\'ik--Besicovitch
theorem are closer to the real, one-dimensional case. Only the
analogue of Jarn\'ik's theorem on badly approximable numbers remains
unchanged and in this case there is a substantial difference in the
proofs of the two theorems.

\section{Applications}
\label{sec:applns}

The connection between the physical phenomenon of resonance and
Diophantine equations can give rise to the notorious problem of small
denominators in which solutions to a variety of questions contain
denominators that can become arbitrarily small.  When these small
denominators are related to very well approximable points, it is
sometimes possible to impose appropriate Diophantine conditions which
overcome the problem by excluding the offending denominators without
significantly affecting the validity of the solution.  The techniques
developed in the metrical theory of Diophantine approximation lend
themselves to this and in particular the Jarník--Besicovitch theorem
allows the determination of the Hausdorff dimension of the associated 
exceptional sets.  Some examples of problems involving small
denominators and the associated exceptional sets are now discussed.
We begin with a very simple example.

\subsection{Partial differential equations}
\label{sec:pdes}

Diophantine approximation has been applied to the wave equation
(\cite{MR50:7828} and more recently \cite{MR95c:35030}), as well as to
the Schrödinger equation \cite{kristensen:_dioph_schroed}. For an
extensive treatment of Diophantine problems related to partial
differential equations, the reader is referred to \cite{Ptashnik84}.

We will illustrate the applications of Diophantine approximation by
Gaussian rationals by considering the following innocuous first-order
linear complex partial differential equation,
\begin{equation}
  \label{eq:6}
  \alpha \dfrac{\partial f(z,t)}{\partial t} + \beta
  \dfrac{\partial f(z,t)}{\partial z} = g(z,t),  
\end{equation}
where $z \in \{x+iy \in \mathbb{C} : x,y \geqslant 0\}, t \geqslant 0$
and $\alpha,\beta$ are non-zero complex numbers. That is, we are
studying the partial differential equation on the interior of the set
defined above under the additional assumption that the functions
involved as well as all their derivatives may be extended to the whole
set. Assume that $g(z,t)$ is smooth (\emph{i.e.}, $C^\infty$) and can
be expressed in the form
\begin{displaymath}
  g(z,t) = \sum_{a,b,c,d \in \mathbb{Z}} g_{a,b,c,d} \exp\left((a+ib)z +
    (c+id)t\right), \quad g_{a,b,d,c} \in \mathbb{C}.
\end{displaymath}
We seek smooth solutions to this equation of the same form, namely
\begin{equation}
  \label{eq:fzt}
  f(z,t) = \sum_{a,b,c,d \in \mathbb{Z}} f_{a,b,c,d}
  \exp\left((a+ib)z + (c+id)t\right), \quad f_{a,b,d,c}
  \in \mathbb{C}.
\end{equation}
Thus, we are not just looking for solutions but rather trying to solve
the partial differential equation subject to boundary conditions. 

As usual, we solve the problem formally by substituting these two
expressions into~\eqref{eq:6} and identifying coefficients on either
side of the equality. Isolating the coefficients of $f$, we get
\begin{equation}
  \label{eq:fabcd}
  f_{a,b,c,d} = \dfrac{\frac{1}{\alpha}}{\frac{\beta}{\alpha}(a+ib)
    + (c+id)} g_{a,b,c,d}. 
\end{equation}
We need the coefficients to decay fast enough so that both the series
\eqref{eq:fzt} and its derivatives are convergent. Since $g$ is
already smooth, the coefficients $g_{a,b,c,d}$ decay rapidly and so
are not obstructing this convergence. However, the denominator of the
fraction may become small and cause the the coefficients $f_{a,b,c,d}$
to become large enough to pose a problem. In order to avoid this, we
see that it is certainly sufficient for the denominator to be bounded
from below by some polynomial in $a,b,c,d$, \emph{i.e.}, we require
for some $K,v > 0$ and for all $(a,b,c,d) \in \mathbb{Z}^4 \setminus
\{0\}$, 
\begin{equation}
  \label{eq:9}
  \abs{\dfrac{\beta}{\alpha}(a+ib) + (c+id)} \geqslant K
  \max\left\{\abs{a},\abs{b}, \abs{c},\abs{d}\right\}^{-v}.
\end{equation}
Since we are only concerned with small denominators, we can assume
without loss of generality that $\abs{a+ib} \asymp \abs{c+id}$, so
that after adjusting $K$ we can drop the dependence on $c$, $d$ on the
right-hand side of \eqref{eq:9} and require
\begin{displaymath}
  \left\vert \dfrac{\beta}{\alpha}-\dfrac{p}{q}\right\vert \geqslant 
  \dfrac{K}{\vert q\vert ^{v+1}},
\end{displaymath}
where $p=a+ib$ and $q=c+id$. Thus we require $\beta/\alpha$ to be of
complex Diophantine type $(K,v)$ for some $K,v$.  The complement of
this set is $ E^*=\cap_{v>1} E_v^*$, where
\begin{equation*}
  E_v^*=\left\{z\in\mathbb{C}\colon  {\text{ for any }} \ K>0, 
    \left\vert z-\dfrac{p}{q}\right\vert  <
    \dfrac{K}{\vert q\vert^{v+1}}  {\text{ for some }} \
    p, q \in\mathbb{Z}[i] \right\}.
\end{equation*}
But it can be readily verified by an argument similar to that giving
the inclusion~\eqref{eq:UEWincl}, that for each $\varepsilon>0$,
\begin{equation*}
  W^*_{v+\varepsilon}\subset E_v^* \subset W^*_v,
\end{equation*}
whence by the properties of Hausdorff dimension given in
\S\ref{sec:HD} and the Jarník--Besicovitch theorem for Gaussian
rational approximation (Theorem~\ref{thm:complexJBt}),
\begin{displaymath}
  \hdim E_v^*=\dfrac{4}{v+1}
\end{displaymath}
for $v\geqslant 1$ and so $\hdim E^*=\lim_{v\to\infty} \hdim E_v^*=0$.
Thus the exceptional set associated with the inequality~\eqref{eq:9}
failing to hold has Hausdorff dimension zero.

\subsection{The rotation number}
\label{sec:rotno}

The rotation number $\rho(f)$ is a measure of how far `on average' a
continuous, orientation preserving homeomorphism $f\colon
\mathbb{S}^1\to \mathbb{S}^1$ moves a point round the circle. We will
not give the fairly lengthy definition which is explained
in~\cite{MDAM,mmdnewton,KatokHass,NiteckiDD} but content ourselves
with the observation that the rotation number of a rotation $r_\alpha$
by an angle $2\pi\alpha$, where $0\leqslant \alpha <1$, given by
\begin{equation*}
  r_\alpha(z)= ze^{2\pi i \alpha}  
\end{equation*}
is, naturally enough, $\alpha$. The rotation number is a nice example
of how Diophantine properties can arise in analysis as it can be shown
that $\rho(f)$ is irrational if and only if $f$ has no periodic points
(see~\cite[Chap.~11,12]{KatokHass} or~\cite[Chap.~1]{NiteckiDD} for
more details).

If $\rho(f)$ is irrational then for $z\in \mathbb{S}^1$, the closure
$A$ of the orbit
\begin{displaymath}
  \omega(z)=\{f^n(z)\colon n\in\mathbb{N}\}
\end{displaymath} 
does not depend on $z$ and either $A$ is perfect and nowhere dense or
$A=\mathbb{S}^1$.  In the latter case $f$ is transitive and is
topologically conjugate to the rotation $r_{\rho(f)}\,$ by the
rotation number $\rho(f)$ of $f$, \emph{i.e.}, there exists an
orientation preserving homeomorphism $\varphi\colon \mathbb{S}^1\to
\mathbb{S}^1$ such that
\begin{equation*}
  f=\varphi^{-1}\circ r_{\rho(f)} \circ \varphi, 
\end{equation*}
usually written $f\sim r_{\rho(f)}$.  This can be regarded as
obtaining a normal form for $f$ and is analogous to 
diagonalising a matrix. 

Denjoy showed that when $f$ is $C^2$ and $\rho(f)$ is irrational, then
$f$ is topologically conjugate to the rotation by $\rho(f)$. More
subtle aspects arise when additional differentiability conditions are
imposed on the conjugation. For example, every $C^\infty$
diffeomorphism $f$ of the the circle is $C^\infty$ conjugate to a
rotation if and only if the rotation number $\rho(f)$ of $f$ is of
Diophantine type~\cite{Yoccoz84}.  In the analytic case, the rotation
numbers of real analytic diffeomorphisms form a set lying strictly
between $\mathcal{D}$ and the set of Bruno
numbers~\cite[p.~92]{Yoccoz95a}.  As in the preceding example, the
Diophantine condition arises from the denominator $1-e^{2\pi i \rho
k}$ in the coefficients for a Fourier series solution of a linearised
auxiliary equation in an iterative Newton's tangent method argument,
modified at each step to retain convergence. In order to guarantee
convergence of the iterative argument and of the Fourier series, the
inequalities
\begin{equation*}
  \vert 1-e^{2\pi i  \rho k}\vert 
  \geqslant 2\left\vert \sin\left(\dfrac{\rho k-j}{2}\right)\right\vert 
  \geqslant \dfrac{2\vert \rho k-j\vert}{\pi},
\end{equation*}
where $\rho k-j\in [0,2\pi)$, must be set against the very rapid decay
of the corresponding Fourier coefficient $f_k$ ($\ll k^{-N}$ for any
$N>0$) in the numerator.  It suffices that $\rho$ is of Diophantine
type $(K,v)$ for some $K>0$, $v>1$, since then
\begin{equation*}
  \left\vert k\rho-j\right\vert \geqslant \dfrac{K}{k^{v}}, \
  \emph{i.e.}, \left\vert \rho-\dfrac{j}{k}\right\vert \geqslant
  \dfrac{K}{k^{v+1}}, 
\end{equation*}
for all $j/k\in\mathbb{Q}$.

Now we saw in~\S\ref{sec:Dt} that when $v>1$, almost all real numbers
are of Diophantine type $(K,v)$ for some positive $K$, \emph{i.e.},
the set
\begin{equation*}
  \mathcal{D}_v = \bigcup_{K>0} \{\alpha\in\mathbb{R}\colon \vert
  \alpha-p/q\vert \geqslant Kq^{-1-v} \text{ for each }
  p/q\in\mathbb{Q}\}  
\end{equation*}
is of full Lebesgue measure. Thus the complementary set $E_v$
\begin{equation*}
  E_{v}= \bigcap_{K>0} \{\alpha\in\mathbb{R}\colon \vert
  \alpha-p/q\vert < Kq^{-2-v} {\text{ for some }} p/q\in\mathbb{Q}\}
\end{equation*}
is null for $v>1$ (see \S\ref{sec:Dt}). As in the preceding section,
its Hausdorff dimension can be determined using the
inclusions~\eqref{eq:UEWincl} and the Jarník--Besicovitch theorem
(Theorem~\ref{thm:JBthm})
\begin{equation*}
  \hdim E_v=\hdim W_{v}=\dfrac{2}{v+1} 
\end{equation*}
for $v\geqslant 1$.  If the rotation number of the smooth circle
function $f$ does not lie in $E=\bigcap_{v>1} E_v=\lim_{v\to\infty}
E_v$ ($E_v$ decreases as $v$ increases), then $f$ is smoothly
conjugate to a rotation. The Hausdorff dimension of the exceptional
set is
\begin{equation}
  \label{eq:dimE0}
  \hdim E= \lim_{v\to\infty} \dfrac{2}{v+1} 
  = 0,
\end{equation}
\emph{i.e.}, the complement of $\mathcal{D}$ has Hausdorff dimension
0.

\subsection{The structure of Julia and Fatou sets}
\label{sec:struct-julia-fatou}

Let $R(z) = P(z)/Q(z)$ be a rational map on the Riemann sphere
$\mathbb{C}_\infty$. A famous result due to Sullivan
\cite{MR99e:58145, MR87i:58103, MR87i:58104} (see also
\cite{MR88e:00019}) states that the Fatou set $F_R = \mathbb{C}_\infty
\setminus J_R$ of such a map has countably many periodic connected
components. These connected components were further classified
according to the type of periodic cycles they are associated with. As
the Julia set is the complementary set of the Fatou set, this
classification also deals with the structure of the Julia set on which
the dynamics of $Q$ is chaotic.

Cycles may be classified according to the value of $dR/dz$ on the
points of the cycle. By the chain rule, this value remains constant,
$\lambda$ say. The cycle is attracting (resp.  repelling) as
$\abs{\lambda} < 1$ (resp. $\abs{\lambda} > 1$). Repelling cycles are
part of the Julia set $J_R$, so no parts of the Fatou set corresponds
to this case. For attracting cycles, an associated periodic connected
component of the Fatou set is in fact the immediate basin of
attraction of this cycle, \emph{i.e.} the union of the connected
components of the Fatou set containing the points of the cycle in
question. 

When $\abs{\lambda}=1$, Diophantine properties of $\lambda$ determine
the behaviour of the dynamics of $R$ and hence the structure of the
Julia and Fatou sets. In this case, $\lambda = \exp(2 \pi i
\alpha)$. When $\alpha$ is rational, the cycle is said to be
parabolic and the associated component of the Fatou set is again the
immediate basin of attraction.

In the case when $\alpha$ is irrational, one is interested in the
situation when the domain of the Fatou set corresponding to the
periodic point is a collection of Siegel discs associated with the
cycle. That is, we are looking for sets on which the dynamics are
topologically conjugate to a rotation of a disc. Such rotation numbers
were studied in the preceding section, where the exceptional set
associated with the failure of this condition was shown to have
Hausdorff dimension zero.

\subsection{Linearising diffeomorphisms}
\label{sec:lindiffeos}
Suppose the complex analytic diffeomorphism $f:\mathbb{C}^{\,n}\to
\mathbb{C}^{\,n}$ has a fixed point.  Without loss of generality, this
can be taken to be the origin, so that $f(0)=0$.  If in a
neighbourhood of $0$, $f$ is analytically (biholomorphically)
conjugate to its linear part or Jacobian $Df\vert_0=A$ say, the function
$f$ is said to be linearisable.  The linearising transformation $\phi$
is given by the solution to the functional equation
\begin{equation*}
  f=\phi^{-1}\circ A  \circ  \phi,
\end{equation*}
known as Schr\"oder's equation when $n=1$.  Thus linearisation is
similar to conjugating a circle map to a rotation, discussed above
in~\S\ref{sec:rotno}.  And problems of small denominators arise when
the eigenvalues $\alpha_1,\dots,\alpha_n$ of $A$ are close to being
resonant in the sense that they are close to satisfying the equation
\begin{equation*}
  \alpha_k=\prod_{r=1}^n \alpha_r^{j_r}
\end{equation*}
for all $\mathbf{j}=(j_1,\dots,j_n)$ with $j_r \in \mathbb{N} \cup
\{0\}$, $r=1,\dots,n$ and $\vert \mathbf{j}\vert_1=\sum_r \vert
j_r\vert  \geqslant 2$. Linearisation is well understood when $n=1$
and the diffeomorphism $f:\mathbb{C}\to \mathbb{C}$ with $f(0)=0$ can
be linearised when $\vert (Df\vert_0)\vert =\vert f'(0)\vert  \not=1$.
The interesting case when $\vert f'(0) \vert =1$ is closely related
via lifts to the conjugacy of a circle map to a rotation, discussed
in~\S\ref{sec:rotno} above, and necessary and sufficient conditions
for the linearisation of a diffeomorphism
$f\colon\mathbb{C}\to\mathbb{C}$ are known~\cite{Yoccoz95a}.  On the
other hand, when $n \geqslant 2$, the problem of finding which 
functions can be linearised is very difficult but Siegel's normal
forms theorem~\cite{Siegel42,Siegel52} gives sufficient conditions on
$Df\vert_0$ for the existence of a linearising transformation $\phi$.  The
point $(\alpha_1,\dots,\alpha_n)$ in $\mathbb{C}^n$ is said to be of
\emph{multiplicative type} $(K,v)$~\cite[p.~191]{ArnoldGM} if
\begin{equation}
  \label{eq:multype}
  \left\vert \alpha_k-\prod_{r=1}^n \alpha_r^{j_r}\right\vert \geqslant 
  K\vert \mathbf{j}\vert_1^{-v} 
\end{equation}
for all $\mathbf{j}\in (\mathbb{N}\cup \{0\})^n$ with
$\vert \mathbf{j}\vert_1\geqslant 2$.  Siegel showed that if the vector
$(\alpha_1,\dots,\alpha_n)$ of eigenvalues of $Df\vert_0$ is of
multiplicative type $(K,v)$ for some $K>0$ and $v>0$, then $f$ can be
linearised in a neighbourhood (further details are
in~\cite{ArnoldGM,Herman87,MoserSRM}).  To stop Siegel's condition
being too restrictive, one chooses $v>(n-1)/2$, since then the set of
points of multiplicative type $(K,v)$ has full measure for any $K>0$.
However, the neighbourhood of linearisation decreases as $v$ increases
(it also depends on $K$ but less significantly) and so we do not want
$v$ to be too large.

Let $\mathcal{E}_v$ denote the exceptional set of points in
$\mathbb{C}^{\,n}$ (regarded as $\mathbb{R}^{2n}$) which for a given
exponent $v$, fail to be of multiplicative type $(K,v)$ for any $K>0$
and so fail to satisfy the conditions of Siegel's theorem and suppose
if $v > (n-1)/2$.  Then $\mathcal{E}_v$ is null and its Hausdorff
dimension is given by~\eqref{eq:dimLv}; namely,
\begin{equation}
  \label{eq:dimEv}
  \dim \mathcal{E}_v = 2(n-1) + \dfrac{n+1}{v+1}.  
\end{equation} 
This result is established by means of an exponential map which
preserves the Hausdorff dimension and allows $\mathcal{E}_v$ to be
replaced by a set involving a more general kind of additive type (see
next section) with a simpler structure~\cite{DRV94}.

\subsection{Lyapunov stability of vector fields}
\label{sec:lypstab}
Consider the differential equation
\begin{equation}
  \label{eq:de}
  \dot{\mathbf{z}}=A\mathbf{z} + Q(\mathbf{z})\in\mathbb{C}^n,
\end{equation}
where $A$ is a $n\times n$ complex matrix and the holomorphic
functions $Q\colon \mathbb{C}^n\to\mathbb{C}^n$ and $\partial
Q_k/\partial z_j$ vanish at the origin $0$.  The obvious solution
$\mathbf{z}_0(t)=0$ is said to be \emph{future (resp. past) stable} if
points near $0$ remain there under evolution by~\eqref{eq:de} to the
future (resp. past).  More precisely, the solution $\mathbf{z}(t)$ is
future (resp. past) stable if for every neighbourhood $N$ of $0$,
there exists a subneighbourhood $N'$ with $0\in N'\subset N$ such that
$\mathbf{z}(0)\in N'$ guarantees that $\mathbf{z}(t)\in N$ for all
$t>0$ (resp. $t<0$).  By a well known theorem of
Lyapunov~\cite{MoserSRM}, the solution is future stable if the real
parts of the eigenvalues of $A$ are at most $0$ and past stable if the
real parts are at least $0$. Thus for the solution $\mathbf{z}(t)$ to
be future and past stable or simply stable, the eigenvalues must have
zero real part and so must be purely imaginary.  The stability of the
solution $\mathbf{z}(t)$ is determined by a remarkable theorem due to
Carath{\'e}odory and Cartan~\cite{MoserSRM} which asserts that
stability is equivalent to $A$ being diagonalisable with purely
imaginary eigenvalues and the vector field being holomorphically
linearisable in a neighbourhood of the origin.   By Siegel's normal
form theorem~\cite{Siegel42,Siegel52}, also used in the preceding
section, this last condition holds if the eigenvalues
$\gamma_k=i\lambda_k$, $\lambda_k\in\mathbb{R}$, satisfy for some
$K$>0, $v >n-1$,
\begin{equation}
  \label{eq:addtype}
  \abs{\lambda_k-\sum_{r=1}^n \lambda_r\, j_r}\geqslant 
  K\abs{\mathbf{j}}_1^{-v}
\end{equation}
for all $\mathbf{j}\in (\mathbb{N}\cup \{0\})^n$ with
$\vert \mathbf{j}\vert_1\geqslant 2$ (note that this is the additive
form of~\eqref{eq:multype}).  The complement of this set of points
$\boldsymbol{\lambda}\in\mathbb{R}^n$ of additive type is the set of
$\boldsymbol{\alpha}\in\mathbb{R}^n$ such that for any $K>0$,
\begin{equation}
  \label{eq:notaddtype}
  \abs{\alpha_k-\sum_{r=1}^n \alpha_r\, j_r}< K\abs{\mathbf{j}}_1^{-v}
\end{equation}
for some $\mathbf{j}\in\mathbb{Z}^n$. This set, $\widehat{E}_v$ say,
is related to and has the same metrical character as the set
\begin{equation}
  \label{eq:hLv}
  \widehat{L}_v=\{\boldsymbol{\alpha}\colon 
  \abs{\mathbf{q}\cdot\boldsymbol{\alpha}}<\abs{\mathbf{q}}_\infty^{-v} 
  \text{ for infinitely many } \ \mathbf{q}\in\mathbb{Z}^n\}
\end{equation}
and in fact by an argument similar to that giving~\eqref{eq:UEWincl},
for any $\varepsilon>0$,
\begin{equation*}
  \widehat{L}_{v+\varepsilon}\subset \widehat{E}_v\subset
  \widehat{L}_v  
\end{equation*}
(see~\cite[Sect.~7.5.2]{MDAM}).  This inclusion implies that the two
sets $ \widehat{L}_v$ and $\widehat{E}_v$ have the same Hausdorff
dimension.  The set $\widehat{L}_v $ is also related to $L_v$ and
roughly speaking, has one degree of freedom less.  The Hausdorff
dimension of $\widehat{L}_v $ is a special case of an `absolute value'
analogue, proved by Dickinson~\cite{HD93}, of the general form of the
Jarník--Besicovitch theorem:
\begin{equation}
  \label{eq:dimhL}
  \dim \widehat{L}_v  = \dim \widehat{E}_v = n-1+ \dfrac{n}{v+1}
\end{equation}
when $v>n-1$ (see also~\cite{DV}); note that $
\widehat{L}_v=\mathbb{R}^n$ otherwise.  Thus the exceptional set of
eigenvalues for which the solution to~\eqref{eq:de} cannot be shown to
be stable has Hausdorff dimension 0.

\subsection{Kolmogorov-Arnol'd-Moser theory}
\label{sec:KAM}

The stability of the solar system is one of the oldest problems in
mechanics~\cite{Moser78}.  It is of course a special case of the $N$
body problem of understanding the motion of $N$ point masses subject
only to gravitational attraction, with all other forces neglected.
When $N=2$, the solution is well known and the periodic solutions in
which the bodies move in an ellipse about their centre of mass persist
forever.  For $N\geqslant 3$, however, the situation is
extraordinarily complicated and is far from being fully understood,
even for solar systems where the mass $m_N$ of the sun is much greater
than the masses of the $n=N-1$ planets.  If, as a first approximation, 
the centre of mass of the system is assumed to coincide with that of
the sun and if the gravitational interactions between the planets and
other effects are neglected, the system decouples into $n$ two-body
problems, in which each planet describes an elliptical orbit around
the sun, with period $T_j$ say and frequency $\omega_j=2\pi/T_j$,
$j=1,\dots,n$.

For each vector $\omega=(\omega_1,\dots,\omega_n)$ of frequencies in
the $n$-dimensional torus $\mathbb{T}^n=\mathbb{S}^1\times \dots\times
\mathbb{S}^1$, the map $\varphi_{\omega} \colon \mathbb{R} \to
\mathbb{T}^n$ given by
\begin{equation*}
  \varphi_{\omega}(t) = \varphi_{\omega}(0)+ t\omega
\end{equation*}
is a quasi-periodic flow on the torus.  

The case $n=1$ corresponds to uniform motion around a circle and so is
periodic.  When the frequencies are all rational, the flow is
periodic.  If the frequencies are not all rational, then by
Kronecker's theorem~\cite{HW}, the flow winds round the torus, densely
filling a subspace of dimension given by the number of rationally
independent frequencies. Thus when the frequencies are independent,
the closure of $\varphi_{\omega}(\mathbb{R})$ is the torus
$\mathbb{T}^n$ and solutions will persist for ever.  Gravitational
interactions between the planets are represented by a small
perturbation of the original Hamiltonian describing the system.
Stability then reduces to the solutions of the perturbed Hamiltonian
system continuing to wind round a perturbed invariant torus.  Of
course this model is idealised and takes no account of the final fate
of the universe.
           
Details of the history of the solution to this problem are in
\cite[Chap.~1]{MoserSRM}.  Siegel's success in overcoming the related
`small denominator' problem in the linearisation of complex
diffeomorphisms (see~\S\ref{sec:lindiffeos}, \S\ref{sec:lypstab}) was
followed by Kolmogorov's conjecture that quasi-periodic solutions for
a perturbed analytic Hamiltonian system not only existed but were
relatively abundant in the sense that they formed a complicated Cantor
type set of positive Lebesgue measure~\cite{Kolmogorov54}.  This was
proved completely in 1962 by Arnol'd~\cite{Arnold63} and independently
Moser proved an analogous result for sufficiently smooth `twist' maps
\cite{Moser62,SiegelMoserCM}.  The results imply that for planets very
much smaller than the sun and for the majority of initial conditions
in which the orbits are close to co-planar circles, distances between
the bodies will remain perpetually bounded, \emph{i.e.}, the planets
will never collide, escape or fall into the sun.  Further details can
be found in \cite{AKNIII} and \cite{KAMtut}.

The differentiability and Diophantine conditions were relaxed
substantially by R\"ussmann in \cite{Russmann83,Russmann91,Russmann01}
(see also \cite[Sect.~6.3]{APIDS}, \cite{NTDS},
\cite[Chap.~1]{MoserSRM}, \cite{Poschel82}).  Another approach is to
use `averaging' methods~\cite{ArnoldGM}; this can involve Diophantine
approximation on manifolds~\cite{DRV89a}, see also~\cite{Russmann01}.

As in the above examples, it turns out that in order to ensure
convergence of a Fourier series and an infinite dimensional extension
of Newton's iterative tangent method, the frequencies
$\omega=(\omega_1,\dots,\omega_n)$ must satisfy a Diophantine
condition, which in this case is
\begin{equation}
  \label{eq:kamda}
  \abs{\mathbf{q}\cdot \boldsymbol{\omega}}=\abs{q_1 \omega_1 + \dots + 
    q_n
    \omega_n} \geqslant K\vert \mathbf{q}\vert_1^{-v}\,, 
\end{equation}
for some positive constants $K=K(\omega)$ and $v=v(\omega)$ for all
non-zero $\mathbf{q}=(q_1,\dots, q_n)$ vectors in $\mathbb{Z}^n$.  The
exponent $v$ is subject to two conflicting requirements.  It should be
large enough ($v>n-1$) to ensure that the Diophantine condition above
is not too restrictive, but small enough to ensure that the
perturbation has physical significance and that the stability is
robust.  The proof breaks down when the frequencies lie in the
complementary exceptional set $E_v$, say, of frequencies which are
close to resonance in the sense that, given any $K>0$, there exists a
$\mathbf{q}\in\mathbb{Z}^n$ such that
\begin{equation*}
  \vert \mathbf{q}\cdot \boldsymbol{\alpha}\vert <K\vert
  \mathbf{q}\vert_1^{-v}.
\end{equation*}
The set $E_v$ is related to the set $\widehat{L}_v$ discussed in the
preceding section and allows us to deduce that
\begin{displaymath}
  \hdim E_v = \hdim \widehat{L}_v=n-1+\dfrac{n}{v+1}
\end{displaymath}
when $v>n-1$. 

Thus Hausdorff dimension plays its part in other branches of
mathematics and mechanics, as well as in number theory. It has even
been of use to Mandelbrot in arousing the interest of mathematicians
in his profoundly original and imaginative
ideas~\cite[Chapter~2]{Cherbit91}.

\section{Acknowledgements}
\label{sec:acknowledgements}

We thank Brent Everitt and Sanju Velani for pointing out the
significance of the Bianchi groups. We also thank the referee
for helpful comments.

\bibliographystyle{amsplain}

\end{document}